\documentclass[francais]{modearticle}

\usepackage[english,frenchb]{babel}
\usepackage[latin1]{inputenc}
 \usepackage[T1]{fontenc}
 \usepackage{lmodern}

% ...

% ...

\equalenv{rem}{rema}

\equalenv{lem}{lemm}
\title
[Th\'eor\`eme de Kaplansky effectif]
{Th\'eor\`eme de Kaplansky effectif pour des valuations de rang 1 centr\'ees sur des anneaux locaux r\'eguliers et complets.}

\author{\firstname{Jean-Christophe}  \lastname{San Saturnino}}

\address{Universit\'e Toulouse III Paul Sabatier\\
Institut de Math\'ematiques de Toulouse\\ 
118 route de Narbonne\\
31062 Toulouse cedex 9 (France)}

\email{san@math.ups-tlse.fr}
\urladdr{http://www.math.univ-toulouse.fr/$\sim$san/}
  

\begin{document}
%% Rsum
\begin{abstract}
On montre que tout anneau local r\'egulier complet muni d'une valuation de rang $1$ peut \^etre plong\'e, en tant qu'anneau valu\'e, dans un anneau de s\'eries de Puiseux g\'en\'eralis\'ees.
\end{abstract}

%% Rsum anglais

\maketitle

\setcounter{section}{-1}
\section{Introduction}
On sait, depuis Newton et Puiseux que, pour un corps $k$ de caract\'eristique $0$ et alg\'ebriquement clos, alors $k(t)$ peut \^etre plong\'e dans le corps $\bigcup\limits_{i\geqslant 1} k\left( \left( t^{1/i}\right) \right) $ \textit{des s\'eries de Puiseux} qui est alg\'ebriquement clos. De plus, si $k$ est muni de la valuation triviale et $k(t)$ de la valuation $t$-adique, on peut alors munir le corps des s\'eries de Puiseux d'une valuation de telle sorte que la restriction \`a $k(t)$ soit la valuation $t$-adique: c'est un exemple d'\textit{extension maximalement compl\`ete} (voir \cite{krull} et \cite{poonen}). 
\\Krull (\cite{krull}) montra, \`a l'aide du Lemme de Zorn, que tout corps muni d'une valuation poss\`ede une extension maximale et que tout corps de s\'eries de Puiseux, muni de sa valuation naturelle, est maximale. L'existence et l'unicit\'e de cette extension maximale fut pos\'ee par Kaplansky (\cite{kaplansky}) qui la d\'emontra en caract\'eristique nulle ainsi que sa non-unicit\'e en caract\'eristique positive. De plus, Poonen (\cite{poonen}) a montr\'e que si le groupe des valeurs de la valuation est divisible et si le corps est alg\'ebriquement clos, alors l'extension maximalement compl\`ete est alg\'ebriquement close. 
\\La question qui vient alors naturellement est: quelle est la forme de cette extension ? En caract\'eristique positive, on sait qu'elle n'est pas de la forme $\bigcup\limits_{i\geqslant 1} k\left( \left( t^{1/i}\right) \right) $ puisque l'\'equation d'Artin-Schreier n'y poss\`ede aucune solution (voir \cite{abhy}, \cite{chevalley}). Il est alors naturel de consid\'erer des anneaux de s\'eries g\'en\'eralis\'ees o\`u les puissances de $t$ varient sur un ensemble bien ordonn\'e. De tels anneaux sont appel\'es des \textit{anneaux de Mal'cev-Neumann} introduits en premier par Hahn en 1908 puis \'etudi\'es par Krull en 1932 (voir \cite{krull}). 
\\ \indent En 1942, Kaplansky (\cite{kaplansky}) montre que tout corps muni d'une valuation ayant un groupe des valeurs divisible et un corps r\'esiduel alg\'ebriquement clos se plonge dans une extension maximalement close. Remarquons que deux cas se pr\'esentent: ou bien la restriction \`a $\mathbb{Q}$ ou $\mathbb{F}_{p}$ est la valuation triviale (cas \'equicaract\'eristique), ou bien la restriction \`a $\mathbb{Q}$ est la valuation $p$-adique (cas mixte). Il a \'egalement montr\'e que dans le cas \'equicaract\'eristique, l'extension maximale est un anneau de Mal'cev-Neumann. 
\\En 1993, Poonen (\cite{poonen}) d\'ecrit explicitement les extensions dans les deux cas. Si $(k,\nu)$ est un corps valu\'e de groupe des valeurs $\Gamma$ divisible et de corps r\'esiduel $k_{\nu}$ alg\'ebriquement clos, alors il existe des plongements dans des anneaux de Mal'cev-Neumann maximalement complets:
\begin{enumerate}
\item $k\hookrightarrow k_{\nu} \left( \left( t^{\Gamma}\right) \right)$ (cas \'equicaract\'eristique);
\item $k\hookrightarrow W(k_{\nu}) \left( \left( p^{\Gamma}\right) \right)$, o\`u $ W(k_{\nu})$ est l'anneau des vecteurs de Witt de $k_{\nu}$ (cas mixte) .
\end{enumerate}
Dans tous les cas les preuves ne construisent pas explicitement le plongement.
\\ \indent Depuis quelques ann\'ees, la th\'eorie des valuations reprend une place importante dans la r\'esolution des singularit\'es et notamment dans l'uniformisation locale des sch\'emas quasi-excellents ,voir par exemple \cite{spiva}. Les questions d'avoir un Th\'eor\`eme de Kaplansky et de conna\^itre explicitement le plongement se posent alors naturellement, on conna\^it l'int\'er\^et d'avoir une param\'etrisation de Puiseux pour l'uniformisation locale des courbes sur un corps de caract\'eristique $0$.
\\ \indent Dans cet article on se propose donc de d\'ecrire de mani\`ere explicite un plongement d'un anneau local r\'egulier et complet muni d'une valuation de rang $1$ \`a l'aide des polyn\^omes-cl\'es d\'efinis dans \cite{spivaherrera} et \cite{mahboub}, r\'esultat qui g\'en\'eralise ceux de Kaplansky et Poonen.
\\ \\ \indent Dans la premi\`ere et la deuxi\`eme partie, nous d\'efinissons les anneaux de Mal'cev-Neumann en suivant \cite{kedlaya} et \cite{poonen} et nous en construisons deux explicitement. 
\\Dans la troisi\`eme partie, nous proposons quelques rappels sur les polyn\^omes-cl\'es d\'efinis pour des valuations de rang $1$. Cet outil est essentiel car il permet de conna\^itre la valuation seulement par la connaissance de la collection bien ordonn\'ee des polyn\^omes-cl\'es, collection qui existe d'apr\'es \cite{spivaherrera} et \cite{spiva}.
\\Dans la quatri\`eme partie, nous \'enon\c{c}ons et d\'emontrons de mani\`ere effective le Th\'eor\`eme de plongement de Kaplansky pour des anneaux locaux, r\'eguliers, complets munis d'une valuation de rang $1$.
\\Dans la derni\`ere partie, nous d\'emontrons des r\'esultats de d\'ependance int\'egrale valables en caract\'eristique mixte qui sont \`a rapprocher de ceux de Spivakovsky (\cite{spiva}) dans le cadre de l'uniformisation locale des sch\'emas quasi-excellents.
\\ \\Je tiens \`a remercier M. Spivakovsky pour son aide pr\'ecieuse, ses conseils avis\'es et la libert\'e qu'il me permet d'avoir dans mes recherches.
\\ \\ \noindent\textbf{Notations.} Soit $(R,\mathfrak{m},k)$ un anneau local, r\'egulier, complet et de dimension $n+ 1$. On note:
 \[p=\left \{ \begin{array}{ccc}  1 & \textup{si} & car(k)=0 \\  car(k)  & \textup{si} & car(k)>0 \end{array} \right.\] 
 Si $R$ est de caract\'eristique mixte, on suppose de plus que $p\notin \mathfrak{m}^{2}$. Par le th\'eor\`eme de Cohen, on peut supposer que:
\[ R=\left \{ \begin{array}{ccc}  k\left[ \left[ u_{1},...,u_{n+1}\right] \right] & \textup{si} & car(R)=car(k) \\  W\left[ \left[ u_{1},...,u_{n}\right] \right]  & \textup{si} & car(R)\neq car(k) \end{array} \right.\]
o\`u $W$ est un anneau complet de valuation discr\`ete de param\`etre r\'egulier $p$ et de corps r\'esiduel $k$. On note:
\[ K_{0}=\left \{ \begin{array}{ccc}  k & \textup{si} & car(R)=car(k) \\  Frac(W)  & \textup{si} & car(R)\neq car(k) \end{array} \right.\]
\[ K_{j}=\left \{ \begin{array}{ccc}  k\left( \left( u_{1},...,u_{j+1}\right) \right) & \textup{si} & car(R)=car(k) \\  W\left( \left( u_{1},...,u_{j}\right) \right)  & \textup{si} & car(R)\neq car(k) \end{array} \right.\]
(on notera parfois $K=K_{n}$). 
\\Soit $\nu$ une valuation de $K$, centr\'ee en $R$, de groupe des valeurs $\Gamma$, telle que $\nu_{\vert K_{n-1}}$ soit de rang $1$. On \'ecrit $(R_{\nu},m_{\nu},k_{\nu})$ son anneau de valuation et on suppose que $k_{\nu}$ est alg\'ebrique sur $k$. On note,  $\Gamma_{1}$ le plus petit sous-groupe isol\'e non-nul de $\Gamma$ et:
\[\Gamma_{\mathbb{Q}}=\Gamma\otimes_{\mathbb{Z}}\mathbb{Q};\]
\[\Gamma'=\bigcup\limits_{i\geqslant 1}\frac{1}{p^{i}}\Gamma.\]
Posons $r$ le plus petit $j$ tel que les $\nu(u_{i_{1}}),...,\nu(u_{i_{j}})$ soient $\mathbb{Z}$-lin\'eairements ind\'ependants si $R$ est \'equicaract\'eristique, ou bien le plus petit $j$ tel que les $\nu(p),\nu(u_{i_{1}}),...,\nu(u_{i_{j}})$ soient $\mathbb{Z}$-lin\'eairements ind\'ependants si $R$ est de caract\'eristique mixte. 
\\On supposera alors, quitte \`a renum\'eroter les variables, que:
\\ $\nu(u_{1}),...,\nu(u_{r})$ sont $\mathbb{Z}$-lin\'eairements ind\'ependantes et $\nu(u_{r+1}),...,\nu(u_{n+1})$ sont $\mathbb{Q}$-combinaisons lin\'eaires de $\nu(u_{1}),...,\nu(u_{r})$ si $R$ est \'equicaract\'eristique;
\\$\nu(p),\nu(u_{1}),...,\nu(u_{r})$ sont $\mathbb{Z}$-lin\'eairements ind\'ependantes et $\nu(u_{r+1}),...,\nu(u_{n})$ sont $\mathbb{Q}$-combinaisons lin\'eaires de $\nu(p),\nu(u_{1}),...,\nu(u_{r})$ si $R$ est de caract\'eristique mixte.
\\On note $\nu_{0}$ la valuation mon\^omiale de $R$ associ\'ee \`a $\mathfrak{m}$ (voir \cite{spiva}, D\'efinition 3.10), c'est-\`a-dire, si $f=\sum\limits_{\alpha}a_{\alpha}u^{\alpha}\in R$ o\`u $\alpha$ est un multi-indice, $a_{\alpha}\in k$ (resp. $a_{\alpha}\in W$) $u^{\alpha}=u_{1}^{\alpha_{1}}...u_{n+1}^{\alpha_{n+1}}$ (resp. $u^{\alpha}=u_{1}^{\alpha_{1}}...u_{n}^{\alpha_{n}}$) et $\mathfrak{m}=(u_{1},...,u_{n+1})$ (resp. $\mathfrak{m}=(p,u_{1},...,u_{n})$), alors:
\[\nu_{0}(f)=\min\left\lbrace \sum\limits_{i=1}^{n+1}\alpha_{i}\nu(u_{i})\:\vert\:a_{\alpha}\neq 0\right\rbrace\] \[\left(\textup{resp.}\: \nu_{0}(f)=\min\left\lbrace \alpha_{0}\nu(p)+\sum\limits_{i=1}^{n}\alpha_{i}\nu(u_{i})\:\vert\:a_{\alpha}\neq 0\right\rbrace\right).\]
Enfin, si $A$ est un anneau, $P,Q\in A\left[ X\right] $ tels que $P=\sum\limits_{i=0}^{n}a_{i}Q^{i}$, $a_{i}\in A[X]$ tels que le degr\'e de $a_{i}$ est strictement inf\'erieur \`a celui de $Q$, on note:
\[d_{Q}^{\:\circ}(P)=n.\]
Si $Q=X$, on notera plus simplement $d^{\:\circ}(P)$ au lieu de $d_{X}^{\:\circ}(P)$.
\\Supposons de plus que $A$ est int\`egre, consid\'erons $w$ une valuation de $A$ et $\gamma\in w(A\setminus\lbrace 0\rbrace)$, on note:\[P'_{\gamma}=\lbrace f\in A\:\vert\:w(f)\geqslant\gamma\rbrace\cup\lbrace 0\rbrace;\]
\[P'_{\gamma,+}=\lbrace f\in A\:\vert\:w(f)>\gamma\rbrace\cup\lbrace 0\rbrace;\]
\[gr_{w}(A)=\bigoplus\limits_{\gamma\in w(A\setminus\lbrace 0\rbrace)}P'_{\gamma}/P'_{\gamma,+};\]
et $in_{w}(f)$ l'image de $f\in A$ dans $gr_{w}(A)$.

\section{Anneaux des s\'eries g\'en\'eralis\'ees} On va d\'efinir des anneaux de s\'eries g\'en\'eralis\'ees (\'egalement appel\'es anneaux de Mal'cev-Neumann) en suivant les constructions donn\'ees par \cite{kedlaya} et \cite{poonen}.
\begin{defi}
Soient $A$ un anneau int\`egre et $G$ un groupe ab\'elien ordonn\'e. On appelle \textbf{anneau des s\'eries formelles g\'en\'eralis\'ees}, not\'e $A\left[ \left[ t^{G}\right] \right] $, l'anneau o\`u les \'el\'ements sont de la forme $\sum\limits_{\gamma\:\in\:G_{+}}a_{\gamma}t^{\gamma}$, avec les  $a_{\gamma}\in A$ tels que l'ensemble $\lbrace \gamma\:\vert\: a_{\gamma}\neq 0\rbrace$ soit bien ordonn\'e.
\\Si $A$ est un anneau int\`egre local de caract\'eristique mixte dont le corps r\'esiduel est de caract\'eristique $p$ et d'id\'eal maximal engendr\'e par $p$, on appelle \textbf{anneau des p-s\'eries formelles g\'en\'eralis\'ees}, not\'e $A\left[ \left[ p^{G}\right] \right] $, l'anneau $A\left[ \left[ t^{G}\right] \right]/N$ o\`u $N$ est l'id\'eal de $A\left[ \left[ t^{G}\right] \right] $ form\'e par les $f=\sum\limits_{\gamma\:\in\:G_{+}}a_{\gamma}t^{\gamma}$ tels que $\sum\limits_{n\in\mathbb{Z}}a_{n+\gamma}p^{n}=0$, pour tout $\gamma\in G$.
\end{defi}

\begin{rem}\textup{ L'anneau $A\left[ \left[ t^{G}\right] \right] $ (resp. l'anneau $A\left[ \left[ p^{G}\right] \right] $) est muni de la \textit{valuation t-adique} (resp. \textit{valuation p-adique})  $v$, \`a valeurs dans $G$, d\'efinie par:
\[ v(f)=inf\lbrace \gamma\:\vert\: a_{\gamma}\neq 0\rbrace\:,\:\forall\: f=\sum\limits_{\gamma\:\in\:G_{+}}a_{\gamma}t^{\gamma}\:\in\:A\left[ \left[ t^{G}\right] \right]  \]
\[(resp.\:\:\forall\: f=\sum\limits_{\gamma\:\in\:G_{+}}a_{\gamma}p^{\gamma}\:\in\:A\left[ \left[ p^{G}\right] \right] ). \]}
\end{rem}

\begin{defi}
Un anneau de s\'eries formelles g\'en\'eralis\'ees (resp. de $p$-s\'eries formelles g\'en\'eralis\'ees) muni de sa valuation $t$-adique (resp. $p$-adique) sera appel\'e un \textbf{anneau de Mal'cev-Neumann}. 
\\Sa valuation $t$-adique (resp. $p$-adique) associ\'ee sera appel\'ee \textbf{valuation de Mal'cev-Neumann}.\label{defmalc}
\end{defi}

\begin{defi}
Soient $f=\sum\limits_{\gamma\:\in\:G_{+}}a_{\gamma}t^{\gamma}\:\in\:A\left[ \left[ t^{G}\right] \right] $ (resp. $ f=\sum\limits_{\gamma\:\in\:G_{+}}a_{\gamma}p^{\gamma}\:\in\:A\left[ \left[ p^{G}\right] \right] $) et $\beta\in G_{+}$, on appelle \textbf{troncature ouverte de f en} $\boldsymbol{\beta}$ la s\'erie g\'en\'eralis\'ee $f(\beta)=\sum\limits_{\gamma<\beta}a_{\gamma}t^{\gamma}$ (resp. $f(\beta)=\sum\limits_{\gamma<\beta}a_{\gamma}p^{\gamma}$) et \textbf{troncature ferm\'ee de f en} $\boldsymbol{\beta}$ la s\'erie g\'en\'eralis\'ee $f[\beta]=\sum\limits_{\gamma\leqslant\beta}a_{\gamma}t^{\gamma}$ (resp. $f[\beta]=\sum\limits_{\gamma\leqslant\beta}a_{\gamma}p^{\gamma}$). Pour $\beta,\beta'\in G_{+}$, $\beta < \beta'$, on note $f[\beta,\beta'[=f(\beta')-f(\beta)$.
\end{defi}

\section{Construction de l'anneau de Mal'cev-Neumann} 
Soient $(R,\mathfrak{m},k)$ un anneau local complet r\'egulier de dimension $n+1$ et $\nu$ une valuation de $K=Frac(R)$, centr\'ee en $R$, de groupe des valeurs $\Gamma$. On va construire un anneau de Mal'cev-Neumann $A_{R}$ dans lequel plonger $R$.
\\\indent Si $R$ est \'equicaract\'eristique, on prend $A_{R}=\overline{k_{\nu}}\left[ \left[ t^{\Gamma'}\right] \right]$, o\`u $\overline{k_{\nu}}$ est une cl\^oture alg\'ebrique de $k_{\nu}$.
\\\indent Si $R$ est de caract\'eristique mixte, on va construire, par r\'ecurrence transfinie, un anneau local $(\overline{W},p\overline{W},\overline{k_{\nu}})$ qui soit une extension de $W$. Dans ce cas, on pose $A_{R}=\overline{W}\left[ \left[ p^{\Gamma'}\right] \right] $.
\\Soit $\overline{k_{\nu}}$ une cl\^oture alg\'ebrique de $k_{\nu}$, on peut la voir comme limite inductive d'extensions alg\'ebriques simples de $k$ puisque $k_{\nu}$ est alg\'ebrique sur $k$. Plus pr\'ecis\'ement, $\overline{k_{\nu}}=k(\left\lbrace \alpha_{i}\right\rbrace_{i\in I} )$ o\`u $I$ est un ensemble bien ordonn\'e et les $\alpha_{i}$ des \'el\'ements alg\'ebriques sur $k$, $i\in I$. Le syst\`eme inductif est alors donn\'e par les inclusions provenant de l'ordre total de $I$. Supposons que $i\in I$ poss\`ede un pr\'ed\'ecesseur imm\'ediat, on est alors emmen\'e \`a consid\'erer une extension de la forme:
\[ \kappa\hookrightarrow\kappa(\alpha) \]
o\`u, par hypoth\`ese de r\'ecurrence, $\alpha$ est alg\'ebrique sur $\kappa$ et $\kappa$ est le corps r\'esiduel d'un anneau local $(A,m_{A})$. Soit $Q$ le polyn\^ome minimal unitaire de $\alpha$ et $P$ un relev\'e unitaire dans $A$. On pose alors $A'=A[\alpha]/(P(\alpha))$ et on a un morphisme d'inclusion:
\[ A\hookrightarrow A' \]
\begin{lem}
$A'$ est un anneau local d'id\'eal maximal $m_{A}A'$ et de corps r\'esiduel $\kappa(\alpha)$.
\end{lem}
\noindent\textit{Preuve}: L'id\'eal $m_{A}A'$ est maximal dans $A'$ car $A'/m_{A}A'\simeq \kappa(\alpha)$. Soit $M$ un autre id\'eal maximal de $A'$, alors $M\cap A=m_{A}$. Pour montrer ceci il suffit de remarquer que $A / (M\cap A)$ est un corps. Soit $a\in A/ (M\cap A)$, $a\neq 0$, l'extension enti\`ere $A\hookrightarrow A'$ induit une extension d'anneaux int\`egres enti\`ere $A/ (M\cap A)\hookrightarrow A'/M$, ainsi $a\in A'/M$ qui est un corps et donc $a^{-1}\in A'/M$. Il existe alors des \'el\'ements $a_{0},...,a_{m-1}\in A/(M\cap A)$ et $m\geqslant 1$ tels que $ a ^{-m}+a_{m-1}a^{-m+1}+...+a_{0}=0$ et donc $a^{-1}=-(a_{m-1}+...+a_{0}a^{m-1})\in A/(M\cap A)$. On remarque enfin que $m_{A}A'=(M\cap A)A'\subset M$ et donc $A'$ est un anneau local.\\ \qed

\noindent Si $i$ est un ordinal limite, notons $\kappa_{l}=k(\left\lbrace \alpha_{j}\right\rbrace_{j\leqslant l} )$, pour tout $l\leqslant i$. On suppose, par hypoth\`ese de r\'ecurrence, que l'on a construit les anneaux locaux $A_{l}$ dont les corps r\'esiduels respectifs sont $\kappa_{l}$ pour tout $l< i$. On pose alors $A_{i}=\bigcup\limits_{l<i}A_{l}$, c'est un anneau local de corps r\'esiduel $\kappa_{i}$. On a donc cr\'e\'e un syst\`eme inductif d'anneaux locaux, on note $\overline{W}$ la limite inductive, c'est un anneau local d'id\'eal maximal $p\overline{W}$, de corps r\'esiduel $\overline{k_{\nu}}$ et on a $W\hookrightarrow \overline{W}$.

\begin{rem}\textup{On a un r\'esultat similaire si $k_{\nu}$ est une extension transcendante de $k$. En effet, si $deg.tr(k_{\nu}\vert k)=l$ alors il existe $t_{1},...,t_{l}$ transcendants tels que $k_{\nu}=k(t_{1},...,t_{l})$. On pose $W'=W[t_{1},...,t_{l}]$ et on consid\`ere l'anneau local $\overline{W'_{pW'}}$, son corps r\'esiduel est  $\overline{k_{\nu}}$ et, $W\hookrightarrow \overline{W'_{pW'}}$.}
\end{rem}

\begin{rem} \textup{Remarquons que $W$ est int\'egralement clos dans $K_{0}$, de plus on a: 
\[W\subsetneqq\overline{W}\left[ p^{\mathbb{Q}}\right] \hookrightarrow\overline{W}\left[ \left[ p^{\Gamma_{\mathbb{Q}}}\right] \right] .\]
Ce dernier morphisme est induit par:
\begin{align*} \mathbb{Q}&\hookrightarrow\Gamma_{\mathbb{Q}}
\\1 &\mapsto\nu(p) \end{align*}}
\end{rem}
\noindent On peut r\'esumer cette section par la proposition suivante:
\begin{prop}
Les anneaux $\overline{k_{\nu}}\left[ \left[ t^{\Gamma'}\right] \right]$ et $\overline{W}\left[ \left[ p^{\Gamma'}\right] \right]$ sont des anneaux de Mal'cev-Neumann au sens de la D\'efinition \ref{defmalc}.
\end{prop}

\section{Rappels sur les polyn\^omes-cl\'es}
On va faire quelques rappels sur les polyn\^omes-cl\'es introduits dans \cite{spivaherrera} et \cite{spiva} pour des valuations de rang $1$. Pour une pr\'esentation plus axiomatique des polyn\^omes-cl\'es, on pourra regarder la pr\'esentation faite par M. Vaqui\'e (\cite{vaquie1}, \cite{vaquie4}, \cite{vaquie2}, \cite{vaquie3}, \cite{vaquie5}). Un lien entre les deux constructions des polyn\^omes-cl\'es est faite dans les travaux de W. Mahboub (\cite{mahboub}).
\\ \indent Soit $K\hookrightarrow K(x)$ une extension de corps simple et transcendante. Soit $\mu'$ une valuation de $K(x)$, notons $\mu:=\mu'_{\vert K}$. On note $G$ le groupe des valeurs de $\mu'$ et $G_{1}$ celui de $\mu$. On suppose de plus que $\mu$ est de rang $1$ et que $\mu'(x)>0$. Enfin, rappelons que pour $\beta\in G$, on note:
\[P'_{\beta}=\lbrace f\in K(x)\:\vert\:\mu'(f)\geqslant\beta\rbrace\cup\lbrace 0\rbrace;\]
\[P'_{\beta,+}=\lbrace f\in K(x)\:\vert\:\mu'(f)>\beta\rbrace\cup\lbrace 0\rbrace;\]
\[gr_{\mu'}(K(x))=\bigoplus\limits_{\beta\in G}P'_{\beta}/P'_{\beta,+};\]
et $in_{\mu'}(f)$ l'image de $f\in K(x)$ dans $gr_{\mu'}(K(x))$.
\begin{defi}
Un \textbf{ensemble complet de polyn\^omes-cl\'es} pour $\mu'$ est une collection bien ordonn\'ee:
\[\textbf{Q}=\lbrace Q_{i}\rbrace_{i\in\Lambda}\subset K[x];\]
telle que, pour tout $\beta\in G$, le groupe additif $P'_{\beta}\cap K[x]$ soit engendr\'e par des produits de la forme $a\prod\limits_{j=1}^{s}Q_{i_{j}}^{\gamma_{j}}$, $a\in K$, tels que $\sum\limits_{j=1}^{s}\gamma_{j}\mu'\left(Q_{i_{j}}\right)+\mu(a)\geqslant\beta$.
\end{defi}
\begin{thm}\textup{(\cite{spivaherrera}, Th\'eor\`eme 62})
Il existe une collection $\textbf{Q}=\lbrace Q_{i}\rbrace_{i\in\Lambda}$ qui soit un ensemble complet de polyn\^omes-cl\'es.
\end{thm}
\begin{rem}
\textup{La preuve consiste \`a construire par r\'ecurrence transfinie l'ensemble de polyn\^omes-cl\'es de type d'ordre au plus $\omega\times\omega$.}
\end{rem}
\begin{defi}
Soit $l\in\Lambda$, un indice $i<l$ est dit \textbf{l-essentiel} s'il existe $n\in\mathbb{N}$ tel que $i+n=l$ ou $i+n<l$ et $d_{Q_{i+n-1}}^{\:\circ}(Q_{i+n})>1$. Dans le cas contraire, on dit que $i$ est \textbf{l-inessentiel}.
\end{defi}
Soit $l\in\Lambda$, on note:
\[\alpha_{i}=d_{Q_{i-1}}^{\:\circ}(Q_{i}),\:\forall\:i\leqslant l;\]
\[\boldsymbol{\alpha_{l+1}}=\lbrace \alpha_{i}\rbrace_{i\leqslant \:l};\]
\[\textbf{Q}_{l+1}=\lbrace Q_{i}\rbrace_{i\leqslant \:l};\]
On utilise \'egalement la notation $\overline{\gamma}_{l+1}=\lbrace \gamma_{i}\rbrace_{i\leqslant\:l}$ o\`u les $\gamma_{i}$ sont tous nuls sauf pour un nombre fini d'entres eux, $\textbf{Q}_{l+1}^{\overline{\gamma}_{l+1}}=\prod\limits_{i\leqslant\: l}Q_{i}^{\gamma_{i}}$.
\\Pour $i<l$, on note:
\[i_{+}=\left \{ \begin{array}{ccc}  i+1 & \textup{si \textit{i} est \textit{l}-essentiel} &  \\  i+\omega  & \textup{sinon} &  \end{array} \right.\] 
\begin{defi}
Un multi-indice $\overline{\gamma}_{l+1}$ est dit \textbf{standard par rapport \`a} $\boldsymbol{\alpha_{l+1}}$ si $0\leqslant \gamma_{i}<\alpha_{i_{+}}$, pour $i\leqslant l$ et si $i$ est $l$-inessentiel, l'ensemble $\lbrace j<i_{+}\:\vert\:j_{+}=i_{+},\:\gamma_{j}\neq 0\rbrace$ est de cardinal au plus $1$.
\\ Un \textbf{mon\^ome l-standard en} $\boldsymbol{Q_{l+1}}$ est un produit de la forme $c_{\overline{\gamma}_{l+1}}\textbf{Q}_{l+1}^{\overline{\gamma}_{l+1}}$, o\`u $c_{\overline{\gamma}_{l+1}}\in K$ et $\overline{\gamma}_{l+1}$ est standard par rapport \`a $\boldsymbol{\alpha_{l+1}}$.
\\ Un \textbf{d\'eveloppement l-standard n'impliquant pas} $\boldsymbol{Q_{l}}$ est une somme finie $\sum\limits_{\beta}S_{\beta}$ de mon\^omes $l$-standards n'impliquant pas $Q_{l}$ o\`u $\beta$ appartient \`a un sous-ensemble fini de $G_{+}$ et $S_{\beta}=\sum\limits_{j} d_{\beta,j}$ est une somme de mon\^omes standards de valuation $\beta$ v\'erifiant $\sum\limits_{j} in_{\mu'}(d_{\beta,j})\neq 0$.
\end{defi}
\begin{defi}
Soient $f\in K[x]$ et $i\leqslant l$, un \textbf{d\'eveloppement i-standard de f} est une expression de la forme:
\[f=\sum\limits_{j=0}^{s_{i}}c_{j,i}Q_{i}^{j},\]
o\`u $c_{j,i}$ est un d\'eveloppement $i$-standard n'impliquant pas $Q_{i}$.
\end{defi}
\begin{rem}
Un tel d\'eveloppement existe, par division Euclidienne et est unique dans le sens o\`u les $c_{j,i}\in K[x]$ sont uniques.
\end{rem}
\begin{defi}
Soient $f\in K[x]$, $i\leqslant l$ et $f=\sum\limits_{j=0}^{s_{i}}c_{j,i}Q_{i}^{j}$ un d\'eveloppement $i$-standard de $f$. On d\'efinit la \textbf{i-troncature de} $\boldsymbol{\mu'}$, not\'ee $\mu_{i}'$, comme \'etant la pseudo-valuation:
\[\mu_{i}'(f)=\min_{0\leqslant j\leqslant s_{i}}\lbrace j\mu'(Q_{i})+\mu'(c_{j,i})\rbrace.\]
\end{defi}
\begin{rem}\label{inegalitetronque}
On peut montrer que c'est en fait une valuation. On a de plus:
\[\forall \:f\in K[x],\: i\in\Lambda, \:\mu_{i}'(f)\leqslant \mu'(f).\]
\end{rem}
On termine en donnant la Proposition 10.1 (Corollaire 50 de \cite{spivaherrera}) et le Corollaire 10.15 de \cite{spiva} que nous utiliserons dans les preuves des Propositions \ref{min} et \ref{taylor}.
\begin{prop}\label{prop10.1}
$\forall\:f\in K[x]$, $\forall\:b\in\mathbb{N}$, 
\[\mu'_{i}(f)-\mu'_{i}\left(\partial_{p^b}f\right)\leqslant\dfrac{p^b}{p^{b_{i}}}\left(\mu'(Q_{i})-\mu'(\partial_{p^{b_{i}}}Q_{i})\right),\]
o\`u $\partial_{p^b}=\dfrac{1}{p^b!}\dfrac{\partial^{p^b}}{\partial x^{p^b}}$ et $b_{i}$ le plus petit $c\in\mathbb{N}^*$ qui maximise $\dfrac{\mu'(Q_{i})-\mu'(\partial_{p^c}Q_{i})}{p^c}$.
\\De plus, il existe $b(i,f)\in\mathbb N^*$ calcul\'e en fonction du d\'eveloppement $i$-standard de $f$, tel que:
\[\mu'_{i}(f)-\mu'_{i}\left(\partial_{p^{b(i,f)}}f\right)=\dfrac{p^{b(i,f)}}{p^{b_{i}}}\left(\mu'(Q_{i})-\mu'(\partial_{p^{b_{i}}}Q_{i})\right).\]
\end{prop}

\begin{prop}
Soit $f=\sum\limits_{j=0}^{s_{i}}c_{j,i}Q_{i}^{j}$ le d\'eveloppement $i$-standard de $f\in K[x]$, on pose:
\[S_{i}=\lbrace j\in\lbrace 0,...,s_{i}\rbrace\:\vert\: j\mu'(Q_{i}) + \mu'(c_{j,i})=\mu'_{i}(f)\rbrace.\]
Soit $j\in S_{i}$, \'ecrivons $j$ sous la forme $j=p^{e}u$, o\`u $p$ ne divise pas $u$ si $car(K_{\mu})=p>0$. 
\\Supposons que $p^{e+1}$ divise $j'$, pour tout $j'\in S_{i}$ tels que $j'<j$. $(*)$
\\Alors:
\[\mu'_{i}(f)=\min_{0\leqslant j\leqslant s_{i}}\lbrace\mu'_{i}\left(\partial_{jb_{i}}f\right)+j\left(\mu'(Q_{i})-\mu'(\partial_{b_{i}}Q_{i})\right) \rbrace \]
et le minimum est atteint pour tout les $j\in S_{i}$ v\'erifiant la condition de divisibilit\'e $(*)$ pr\'ec\'edente.\label{coro10.15}
\end{prop}

On va utiliser les polyn\^omes-cl\'es dans le cadre des anneaux locaux r\'eguliers, ils interviennent de mani\`ere fondamentale dans la d\'emonstration du Th\'eor\`eme \ref{thmequi}.
\\

\noindent\textbf{Polyn\^omes-cl\'es dans une tour d'extensions de corps.} Pour $j\in \lbrace r+1,...,n\rbrace$, on note $\lbrace Q_{j,i}\rbrace_{i\in\Lambda_{j}}$ l'ensemble des polyn\^omes-cl\'es de l'extension $K_{j-1}\hookrightarrow K_{j-1}(u_{j})$, $\textbf{Q}_{j,i}=\left\lbrace Q_{j,i'}\vert i'\in\Lambda_{j},i'<i\right\rbrace $, $\Gamma^{(j)}$ le groupe des valeurs de $\nu_{\vert K_{j}}$ et $\nu_{j,i}$ la $i$-troncature de $\nu$ pour cette extension. Soient $\beta_{j,i}=\nu(Q_{j,i})$ et $b_{j,i}$ le plus petit \'el\'ement $b$ de $\mathbb{N}$ qui maximise $\dfrac{\beta_{j,i}-\nu(\partial_{j,p^{b}}Q_{j,i})}{p^{b}}$, o\`u $\partial_{j,s}=\dfrac{1}{s!}\dfrac{\partial^{s}}{\partial u_{j}^{s}}$, $s\in\mathbb N$. Soit $\varepsilon_{j,i}=\dfrac{\beta_{j,i}-\nu(\partial_{j,p^{b_{j,i}}}Q_{j,i})}{p^{b_{j,i}}}$, on a:
\begin{lem}\label{suiteepsiloncroiss}(\cite{spiva}, Lemme 10.4) 
La suite $(\varepsilon_{j,i})_{i}$ est strictement croissante.
\end{lem}
\noindent\textit{Preuve}: Nous reprennons la preuve du Lemme 10.4 de \cite{spiva}. Fixons $j\in\lbrace r+1,...,n\rbrace$ et consid\'erons $i_1$, $i_2\in \Lambda_j$ deux ordinaux. Il faut montrer que:
\[\dfrac{\beta_{j,i_1}-\nu(\partial_{j,p^{b_{j,i_1}}}Q_{j,i_1})}{p^{b_{j,i_1}}}<\dfrac{\beta_{j,i_2}-\nu(\partial_{j,p^{b_{j,i_2}}}Q_{j,i_2})}{p^{b_{j,i_2}}}.\]
On peut supposer que $i_2=i_{1\:+}$ (c'est-\`a-dire $i_2=i_1+1$ ou $i_2=i_1+\omega$), on conclut dans le cas g\'en\'eral par r\'ecurrence transfinie sur $i_2-i_1$. Par la Proposition \ref{prop10.1}, il existe $b(i_1,Q_{j,i_2})$ calcul\'e en fonction du d\'eveloppement $(j,i_1)$-standard de $Q_{j,i_2}$, tel que:
\[\nu_{j,i_1}(Q_{j,i_2})-\nu_{j,i_1}\left(\partial_{p^{b(i_1,Q_{j,i_2})}}Q_{j,i_2}\right)=\dfrac{p^{b(i_1,Q_{j,i_2})}}{p^{b_{i_1}}}\left(\nu(Q_{j,i_1})-\nu(\partial_{p^{b_{i_1}}}Q_{j,i_1})\right).\]
Vu que $d^{\:\circ}\left(\partial_{p^{b(i_1,Q_{j,i_2})}}Q_{j,i_2}\right)<d^{\:\circ}\left(Q_{j,i_2}\right)$, on montre facilement que:
\[\nu_{j,i_1}\left(\partial_{p^{b(i_1,Q_{j,i_2})}}Q_{j,i_2}\right)=\nu\left(\partial_{p^{b(i_1,Q_{j,i_2})}}Q_{j,i_2}\right).\]
Par d\'efinition du d\'eveloppement $(j,i_1)$-standard, on a:
\[\beta_{j,i_2}>\alpha_{j,i_2}\beta_{j,i_1}=\nu_{j,i_1}(Q_{j,i_2}).\]
Ainsi, par d\'efinition de $\varepsilon_{j,i_2}$:
\begin{align*}
\varepsilon_{j,i_1}&=\dfrac{\nu_{j,i_1}(Q_{j,i_2})-\nu_{j,i_1}\left(\partial_{p^{b(i_1,Q_{j,i_2})}}Q_{j,i_2}\right)}{p^{b(i_1,Q_{j,i_2})}}\\&=\dfrac{\alpha_{j,i_2}\beta_{j,i_1}-\nu\left(\partial_{p^{b(i_1,Q_{j,i_2})}}Q_{j,i_2}\right)}{p^{b(i_1,Q_{j,i_2})}}\\&<\dfrac{\beta_{j,i_2}-\nu\left(\partial_{p^{b(i_1,Q_{j,i_2})}}Q_{j,i_2}\right)}{p^{b(i_1,Q_{j,i_2})}}\\&\leqslant\varepsilon_{j,i_2}.
\end{align*}
On en conclut que la suite $(\varepsilon_{j,i})_{i}$ est strictement croissante, pour tout $j\in\lbrace r+1,...,n\rbrace$.
\\\qed

\section{Le th\'eor\`eme de plongement de Kaplansky}
Dans cette section, on suppose que $(R,\mathfrak{m},k)$ est un anneau local complet r\'egulier de dimension $n+1$. Si $R$ est de caract\'eristique mixte, on suppose de plus que $p\notin\mathfrak{m}^{2}$. Notons $\nu$ une valuation de $K=Frac(R)$, centr\'ee en $R$, de groupe des valeurs $\Gamma$ et telle que $\nu_{\vert K_{n-1}}$ soit de rang $1$.

\begin{thm}
Il existe un anneau de Mal'cev-Neumann $A_{R}$ et un monomorphisme d'anneaux \[\iota: R\hookrightarrow A_{R} \] tel que $\nu$ soit la restriction \`a $R$ de la valuation de Mal'cev-Neumann associ\'ee \`a $A_{R}$.
\\Pour $f\in R$, on appelle $\iota(f)$ un \textbf{d\'eveloppement de Puiseux de $f$ par rapport \`a} $\boldsymbol{\nu}$.\label{thmequi}
\end{thm}

\begin{rem}
$ A_{R}=\left \{ \begin{array}{ccc} \overline{k_{\nu}}\left[ \left[ t^{\Gamma'}\right] \right]  & \textup{si} & car(R)=car(k) \\  \overline{W}\left[ \left[ p^{\Gamma'}\right] \right]  & \textup{si} & car(R)\neq car(k) \end{array} \right.$ 
\end{rem}

\begin{rem}
On sait que:
\[ R=\left \{ \begin{array}{ccc}  k\left[ \left[ u_{1},...,u_{n+1}\right] \right] & \textup{si} & car(R)=car(k) \\  W\left[ \left[ u_{1},...,u_{n}\right] \right]  & \textup{si} & car(R)\neq car(k) \end{array} \right.\]
la preuve consiste donc \`a d\'efinir, par r\'ecurrence transfinie, le d\'eveloppement de Puiseux de $u_{1},...,u_{n+1}$ (resp. $p,u_{1},...,u_{n}$) \`a l'aide des polyn\^omes-cl\'es.
\end{rem}

\noindent\textit{Preuve}: On va faire la preuve de ce th\'eor\`eme seulement dans le cas o\`u $R$ est de caract\'eristique mixte. Le cas o\`u $R$ est \'equicaract\'eristique se traite de la m\^eme mani\`ere en rempla\c{c}ant $p$ par $t$ et en prenant les coefficients directement dans $\overline{k_{\nu}}$.

Dans ce qui suit on va construire un d\'eveloppement de Puiseux en lien avec les polyn\^omes-cl\'es. Remarquons que d\'efinir un d\'eveloppement de Puiseux pour un \'el\'ement de $R$ revient \`a d\'efinir $n+1$ s\'eries $\iota(p),\iota(u_{1}),...,\iota(u_{n})$ formellement ind\'ependantes sur $W$.
\\On va construire le morphisme $\iota$ par r\'ecurrence sur $n-r$. Si $n=r$, on pose $\iota(p)=p^{\nu(p)}$ et $\iota (u_{j})=p^{\nu (u_{j})}$, $j\in \lbrace 1,...,r\rbrace$ (remarquons que, pour que $\iota$ soit un morphisme, comme $\nu_{\vert K_{n-1}}$ est de rang $1$, on choisit une fois pour toute un plongement $\Gamma_1\hookrightarrow\mathbb R$ qui envoie $\nu(p)$ sur $1$).
\\Supposons que $n>r$ et que l'on a d\'ej\`a construit un monomorphisme d'anneaux valu\'es:
\[ \iota_{n-1}:R_{n-1}\hookrightarrow \overline{W}\left[ \left[ p^{\Gamma'}\right] \right]  \]
tel que $\nu_{\vert R_{n-1}}$ soit induite par la valuation $p$-adique, $R_{n-1}$ d\'esignant l'anneau $W\left[ \left[ u_1,...,u_{n-1}\right] \right] $. Pour $j\in\lbrace 1,...,n-1\rbrace$, on note $u_{j}(p)=\iota_{n-1}(u_{j})$.
\\Nous allons construire la s\'erie g\'en\'eralis\'ee $u_{n}(p)$ par r\'ecurrence transfinie sur un sous-ensemble bien ordonn\'e de $\Gamma'$.
\\Soit $\beta\in\Gamma_{+}$, on note $i_{\beta}=\min\lbrace i\in\Lambda_{n}\:\vert\:\beta\leqslant\varepsilon_{n,i}\rbrace$ et par convention, si $\lbrace i\in\Lambda_{n}\:\vert\:\beta\leqslant\varepsilon_{n,i}\rbrace=\emptyset$, on prendra $i_{\beta}=\Lambda_{n}$ (c'est-\`a-dire le plus petit ordinal strictement plus grand que n'importe quel \'el\'ement de $\Lambda_{n}$).
\\Supposons donn\'ee une s\'erie g\'en\'eralis\'ee $u_{n}(\beta)=\sum\limits_{\gamma<\beta}a_{\gamma}p^{\gamma}\in\overline{W}\left[ \left[ p^{\Gamma'}\right] \right] $, on consid\`ere le morphisme d'anneaux $\iota_{\beta}:R\rightarrow\overline{W}\left[ \left[ p^{\Gamma'}\right] \right] $ d\'efini par:
\begin{align*}
\iota_{\beta}(p)&=p^{\nu(p)}; \\
\iota_{\beta}(u_{j})&=u_{j}(p),\:\forall\:j\in\lbrace 1,...,n-1\rbrace ;\\
\iota_{\beta}(u_{n})&=u_{n}(\beta).
\end{align*}
\begin{defi}\label{defdevpart}
On dit que $u_{n}(\beta)$ est un \textbf{d\'eveloppement de Puiseux partiel} de $u_{n}$ si les deux conditions suivantes sont v\'erifi\'ees pour tout $j\in\lbrace 1,...,n\rbrace$:
\begin{enumerate}
\item $v(\iota_{\beta}(Q_{n,i}))=\beta_{n,i}$, $\forall i \in \Lambda_{n}$ tel que $i<i_{\beta}$;
\item $v(\iota_{\beta}(Q_{n,i_{\beta}}))\geqslant \min\limits_{q\in\mathbb{N}}\lbrace \nu(\partial_{n,p^{q}}Q_{n,i_{\beta}})+p^{q}\beta\rbrace$ (si $i_{\beta}=\Lambda_{n}$, on consid\`ere cette condition toujours v\'erifi\'ee).
\end{enumerate}
\end{defi}

\indent Soit $T$ une nouvelle variable et consid\'erons le morphisme d'anneaux \\$\iota_{\beta,T}:R\rightarrow\overline{W}\left[ \left[ p^{\Gamma'},T\right] \right] $ d\'efini par:
\begin{align*}
 \iota_{\beta,T}(p)&=p^{\nu(p)}; \\
 \iota_{\beta,T}(u_{j})&=u_{j}(p),\:\forall\:j\in\lbrace 1,...,n-1\rbrace ; \\
 \iota_{\beta,T}(u_{n})&=u_{n}(\beta)+T.
 \end{align*} 
On note $\nu_{\beta}$ l'extension \`a $\overline{W}\left[ \left[ p^{\Gamma'},T\right] \right] $ de la valuation $p$-adique $v$ de $\overline{W}\left[ \left[ p^{\Gamma'}\right] \right] $ telle que $\nu_{\beta}(T)=\beta$ et on suppose que $in_{\nu_{\beta}}(T)$ est transcendant sur $gr_{v}\left( \overline{W}\left[ \left[  p^{\Gamma'}\right] \right] \right) $. On pose alors $\mu_{\beta}=\nu_{\beta\vert R}$ o\`u $R$ est vu comme sous-anneau de $\overline{W}\left[ \left[ p^{\Gamma'},T\right] \right] $ via le monomorphisme $\iota_{\beta,T}$.
\\Supposons que $i_{\beta}=\Lambda_{n}$, alors $\mu_{\beta}=\nu$. Sinon, supposons que $i_{\beta}<\Lambda_{n}$, c'est-\`a-dire qu'il existe $i\in\Lambda_{n}$ tel que $\varepsilon_{n,i}\geqslant\beta$. On note alors:
\[ \Gamma_{\beta}=\Gamma^{(n-1)}\otimes_{\mathbb{Z}}\mathbb{Q} + \sum\limits_{i<i_{\beta}}\mathbb{Q}\beta_{n,i}\subset\Gamma'. \]
\begin{lem}\label{lemmebeta}
On a les assertions suivantes:\begin{enumerate}
\item La valuation $\nu_{\beta\vert K_{n-1}}$ est l'unique valuation telle que:
\[\nu_{\beta\vert K_{n-1}[\textbf{Q}_{n,i_{\beta}}]}=\nu_{\vert K_{n-1}[\textbf{Q}_{n,i_{\beta}}]};\]
\[\nu_{\beta}(Q_{n,i_{\beta}})=\min_{q\in\mathbb{N}}\left\lbrace \nu\left( \partial_{n,p^{q}}Q_{n,i_{\beta}}\right) +p^{q}\beta\right\rbrace ,\]
et $in_{\nu_{{\beta}}}(Q_{n,i_{\beta}})$ est transcendant sur $gr_{\nu_{\beta}}\left( K_{n-1}[\textbf{Q}_{n,i_{\beta}}]\right)$.
\item Consid\'erons les sous-alg\`ebres gradu\'ees $gr_{\mu_{\beta}}\left( K_{n-1}[\textbf{Q}_{n,i_{\beta}}]\right)\subset gr_{\mu_{\beta}}(R)$ et $gr_{\nu_{n,i_{\beta}}}\left( K_{n-1}[\textbf{Q}_{n,i_{\beta}}]\right)\subset gr_{\nu_{n,i_{\beta}}}\left( R\right)$. On a alors un isomorphisme d'alg\`ebres gradu\'ees:
\[gr_{\mu_{\beta}}\left( K_{n-1}[\textbf{Q}_{n,i_{\beta}}]\right)\simeq gr_{\nu_{n,i_{\beta}}}\left( K_{n-1}[\textbf{Q}_{n,i_{\beta}}]\right)\]
qui peut \^etre \'etendu en un isomorphisme entre $gr_{\mu_{\beta}}(R)$ et $gr_{\nu_{n,i_{\beta}}}\left( R\right)$ en envoyant $in_{\mu_{\beta}}(Q_{n,i_{\beta}})$ sur $in_{\nu_{n,i_{\beta}}}(Q_{n,i_{\beta}})$, mais la graduation n'est, en g\'en\'eral, pas pr\'eserv\'ee, sauf si l'une des deux conditions \'equivalentes de (3) est v\'erifi\'ee.
\item $\mu_{\beta}=\nu_{n,i_{\beta}}$ ssi $\beta=\varepsilon_{n,i_{\beta}}$.
\item $\forall\:h\in R$, $\nu_{n,i_{\beta}}(h)\leqslant\nu(h)$.
\\Supposons que $\beta=\varepsilon_{n,i_{\beta}}$ (donc $\mu_{\beta}=\nu_{n,i_{\beta}}$). On a alors, pour tout $h\in R$:
\begin{align*} &\nu_{n,i_{\beta}}(h)=\nu(h) \Leftrightarrow in_{\nu_{n,i_{\beta}}}(h)\notin \ker\left( gr_{\nu_{n,i_{\beta}}}\left( R\right)\rightarrow gr_{\nu}(R)\right)  \\ &\Leftrightarrow in_{\nu_{\beta}}\left( \iota_{\beta,T}(h)\right)\notin \ker \Big( gr_{\nu_{\beta}}\left( \overline{W} \left[ \left[  p^{\Gamma'},T\right] \right]\right)  \rightarrow gr_{v} \left(  \overline{W} \left[ \left[  p^{\Gamma'}\right] \right]\right)   \Big) .     \end{align*}
En particulier, il y a \'egalit\'e si $in_{\nu_{\beta}}(T)$ n'appara\^it pas dans $in_{\nu_{\beta}}( \iota_{\beta,T}(h)) $.
\end{enumerate}
\label{rmqval}\end{lem}
\noindent\textit{Preuve}: (1): $\nu_{\beta\vert K_{n-1}}=\nu_{\vert K_{n-1}}$ par d\'efinition de $\iota_{\beta,T}$, $\nu_{\beta}$ et $v$. Pour $i<i_\beta$, alors, $\beta_{n,i}<\beta$ et comme $u_n(\beta)$ est un d\'eveloppement de Puiseux partiel, on obtient l'\'egalit\'e:
\[\nu_\beta\left(\iota_{\beta,T}(Q_{n,i})\right)=v\left(\iota_{\beta}(Q_{n,i})\right)=\beta_{n,i}.\]
Enfin, comme:
\[Q_{n,i_\beta}\left(u_1,...,u_{n-1},u_n+T\right)=\sum\limits_{l=0}^{d_{u_n}^{\:\circ}(Q_{n,i_\beta})}\partial_{n,l}Q_{n,i_\beta}(u_1,...,u_n)T^l,\]
on en d\'eduit que:
\[\nu_{\beta}(Q_{n,i_{\beta}})=\min_{l\in\mathbb{N}}\left\lbrace \nu\left( \partial_{n,l}Q_{n,i_{\beta}}\right) +l\nu_\beta(T)\right\rbrace ,\]
o\`u le minimum est atteint avec $l=1$ si $car(k)=0$, une puissance de $p=car(k)$ sinon. Ainsi, $in_{\nu_{\beta}}(T)$ appara\^it dans $in_{\nu_{{\beta}}}(Q_{n,i_{\beta}})$. Comme $in_{\nu_{\beta}}(T)$ est transcendant sur $gr_{v}\left( \overline{W}\left[ \left[  p^{\Gamma'}\right] \right] \right) $, on en d\'eduit que $in_{\nu_{{\beta}}}(Q_{n,i_{\beta}})$ est transcendant sur $gr_{\nu_{\beta}}\left( K_{n-1}[\textbf{Q}_{n,i_{\beta}}]\right)$. L'unicit\'e de $\nu_\beta$ v\'erifiant les propri\'et\'es pr\'ec\'edemment d\'emontr\'ees provient de la d\'efinition m\^eme de cette valuation.
\\(2): Par d\'efinition et construction des polyn\^omes-cl\'es et de la valuation tronqu\'ee $\nu_{n,i_\beta}$, on obtient l'\'egalit\'e:
\[\nu_{n,i_\beta\vert K_{n-1}[\textbf{Q}_{n,i_{\beta}}]}=\nu_{\vert K_{n-1}[\textbf{Q}_{n,i_{\beta}}]}\]
qui nous d\'efinit un isomorphisme naturel d'alg\`ebres gradu\'ees:
\[gr_{\nu_{n,i_{\beta}}}\left( K_{n-1}[\textbf{Q}_{n,i_{\beta}}]\right)\overset{\sim}{\longrightarrow}gr_{\mu_{\beta}}\left( K_{n-1}[\textbf{Q}_{n,i_{\beta}}]\right)\]
En envoyant $in_{\nu_{n,i_{\beta}}}(Q_{n,i_{\beta}})$ sur $in_{\mu_{\beta}}(Q_{n,i_{\beta}})$, on prolonge l'isomorphisme pr\'ec\'edent en un isomorphisme entre $gr_{\nu_{n,i_{\beta}}}\left( R\right)$ et $gr_{\mu_{\beta}}(R)$, la graduation \'etant pr\'eserv\'ee seulement si $\mu_\beta=\nu_{n,i_\beta}$.
\\(3): Supposons que $\beta=\varepsilon_{n,i_\beta}$, par d\'efinition des polyn\^omes-cl\'es et de $\mu_\beta$, il suffit de montrer que $\mu_\beta(Q_{n,i_\beta})=\nu_{n,i_\beta}(Q_{n,i_\beta})$. Soit $q_0\in\mathbb N$ tel que $\mu_\beta(Q_{n,i_\beta})=\nu\left( \partial_{n,p^{q_0}}Q_{n,i_{\beta}}\right) +p^{q_0}\beta$. Par d\'efinition de $\mu_\beta$, on a:
\[\mu_\beta(Q_{n,i_\beta})=\nu\left( \partial_{n,p^{q_0}}Q_{n,i_{\beta}}\right) +p^{q_0}\beta\leqslant \nu\left( \partial_{n,p^{b_{n,i_\beta}}}Q_{n,i_{\beta}}\right) +p^{b_{n,i_\beta}}\varepsilon_{n,i_\beta}=\beta_{n,i_\beta}.\]
Par d\'efinition de $\varepsilon_{n,i_\beta}$, on a:
\[\varepsilon_{n,i_\beta}\geqslant\dfrac{\beta_{n,i_\beta}-\nu\left( \partial_{n,p^{q_0}}Q_{n,i_{\beta}}\right)}{p^{q_0}},\]
c'est-\`a-dire:
\[\mu_\beta(Q_{n,i_\beta})=\nu\left( \partial_{n,p^{q_0}}Q_{n,i_{\beta}}\right) +p^{q_0}\beta\geqslant\beta_{n,i_\beta}.\]
R\'eciproquement, si $\mu_\beta=\nu_{n,i_\beta}$, alors:
\[\beta_{n,i_\beta}=\nu\left( \partial_{n,p^{q_0}}Q_{n,i_{\beta}}\right) +p^{q_0}\beta\leqslant\nu\left( \partial_{n,p^{b_{n,i_\beta}}}Q_{n,i_{\beta}}\right) +p^{b_{n,i_\beta}}\beta,\]
ce qui donne:
\[\varepsilon_{n,i_\beta}=\dfrac{\beta_{n,i_\beta}-\nu\left( \partial_{n,p^{b_{n,i_\beta}}}Q_{n,i_{\beta}}\right)}{p^{b_{n,i_\beta}}}\leqslant\beta.\]
Enfin, rappelons que, par d\'efinition de $i_\beta$, $\beta\leqslant\varepsilon_{n,i_\beta}$.
\\(4): Par la Remarque \ref{inegalitetronque}, pour tout $h\in R$, $\nu_{n,i_{\beta}}(h)\leqslant\nu(h)$. La premi\`ere \'equivalence est \'evidente, la deuxi\`eme provient du fait que l'on a suppos\'e $\mu_\beta=\nu_{n,i_\beta}$ et que $\mu_\beta(h)=\nu_\beta(\iota_{\beta,T}(h))$ ainsi que $v(\iota(h))=\nu(h)$.\\\qed

\indent Commen\c{c}ons notre r\'ecurrence transfinie par $\beta=\nu(u_{n})=\beta_{n,1}$. On pose alors $u_{n}(\beta)=0$ et on a $i_{\beta}=1$, $\mu_{\beta}=\nu_{n,i_{\beta}}=\nu_{n,1}$; $u_n(\beta)$ est ainsi un d\'eveloppement de Puiseux partiel de $u_{n}$.
\\Supposons $u_{n}(\beta)$ construit pour un certain $\beta\in\Gamma_{+}$ tel que $\beta > \nu(u_{n})$ et d\'efinissons le coefficient $a_{n,\beta}$ de $p^{\beta}$ de $u_{n}(p)$. On suppose \'egalement, par hypoth\`ese de r\'ecurrence, que $\beta=\varepsilon_{n,i_{\beta}}$ ou que $\beta\in\Gamma_{\beta}$.
\\Si $\beta\notin\Gamma_{\beta}$, comme $\beta=\varepsilon_{n,i_{\beta}}$, alors $\beta_{n,i_{\beta}}\notin\Gamma_{\beta}$ et donc $i_{\beta}=\max\Lambda_{n}$. Dans ce cas on a $\nu=\nu_{n,i_{\beta}}=\mu_{\beta}$ et on pose alors:
\[ u_{n}(p)=u_{n}(\beta)+p^{\beta}. \]
\\Si $\beta\in\Gamma_{\beta}$ alors $\nu_{\beta}(Q_{n,i_{\beta}})=\min\limits_{q\in\mathbb{N}}\lbrace \nu(\partial_{n,p^{q}}Q_{n,i_{\beta}})+p^{q}\beta\rbrace\in\Gamma_{\beta}$ et donc:
\[ \exists\: d\in K_{n-1},l_{1},...,l_{t}\in\Lambda_{n},\lambda\in\mathbb{N},\lambda_{1},...,\lambda_{t}\in\mathbb{Z}\]tels que:\[\lambda\nu_{\beta}(Q_{n,i_{\beta}})=\sum\limits_{j=1}^{t}\lambda_{j}\beta_{n,l_{j}}+\nu(d). \]
On pose alors:
\[ z=\left \{ \begin{array}{ccc}  \dfrac{Q_{n,i_{\beta}}^{\lambda}}{d\prod\limits_{j=1}^{t}Q_{n,l_{j}}^{\lambda_{j}}}\:\textup{mod}\: m_{\nu}\in k_{\nu} & \textup{si} & \beta=\varepsilon_{n,i_{\beta}}\\  0   & \textup{si} & \beta<\varepsilon_{n,i_{\beta}} \end{array} \right.\]
Notons $W_{r+1},...,W_{n-1},W_{n}$ les supports respectifs de $u_{r+1}(p),...,u_{n-1}(p)$, $u_{n}(\beta)$, on note alors $u_{j}(p)=\sum\limits_{\gamma\in W_{j}}a_{j,\gamma}p^{\gamma}$ pour $j\in\lbrace r+1,...,n-1\rbrace$ et $u_{n}(\beta)=\sum\limits_{\gamma\in W_{n}}a_{n,\gamma}p^{\gamma}$.
Posons, pour $j\in\lbrace r+1,...,n\rbrace$, $\mathfrak{a}_{j}=\lbrace a_{j,\gamma}\:\vert\:\gamma\in W_{j}\rbrace\subset\overline{W}$, $\overline{\mathfrak{a}_{j}}=\lbrace \overline{a_{j,\gamma}}\:\vert\:\gamma\in W_{j}\rbrace\subset\overline{k_{\nu}}$ son image modulo $p$, $\mathfrak{a}=(\mathfrak{a}_{r+1},...,\mathfrak{a}_{n})$ et $\overline{\mathfrak{a}}=(\overline{\mathfrak{a}_{r+1}},...,\overline{\mathfrak{a}_{n}})$. Soit $X$ une variable ind\'ependante, si on remplace $T$ par $Xp^{\beta}$ dans $in_{\nu_{\beta}}(\iota_{\beta,T}(Q_{n,i_{\beta}}))$, on obtient un r\'esultat de la forme:
\[ fp^{\nu_{\beta}(Q_{n,i_{\beta}})},\:f\in K_{0}\left[ \mathfrak{a},X\right] \ . \]
On note alors $\overline{f}\in k_{\nu}[\overline{\mathfrak{a}},X]\subset\overline{k_{\nu}}[X]$ l'image de $f$ modulo $p$.
\\De plus, pour $j\in\lbrace 1,...,t\rbrace$, $in_{\nu_{\beta}}(\iota_{\beta}(Q_{n,l_{j}}))$ est de la forme:
\[ c_{j}p^{\beta_ {n,l_{j}}},\:c_{j}\in \overline{W}\]
et $in_{\nu_{\beta}}(\iota_{\beta}(d))$ est de la forme:
\[ \delta p^{\nu(d)},\:\delta\in \overline{W}.\]
Notons $\overline{c_{j}}$ et $\overline{\delta}$ dans $\overline{k_{\nu}}$ les images respectives de $c_{j}$ et de $\delta$ modulo $p$.
\\Ainsi, $\dfrac{\overline{f}^{\lambda}}{\overline{\delta}\prod\limits_{j=1}^{t}\overline{c_{j}}^{\lambda_{j}}}=z$ induit une \'equation alg\'ebrique en $X$ sur $\overline{k_{\nu}}$, on note alors $\alpha_{n,\beta}\in\overline{k_{\nu}}$ une de ses racines et $a_{n,\beta}\in\overline{W}$ un relev\'e. Deux cas se pr\'esentent:
\begin{enumerate}
\item $a_{n,\beta}$ est transcendant sur $K_{0}[\mathfrak{a}]$.
On pose alors \[u_{n}(p)=u_{n}(\beta)+a_{n,\beta}p^{\beta},\] on a $\nu=\nu_{n,i_{\beta}}$ et on arr\^ete l'algorithme.
\item $a_{n,\beta}$ est alg\'ebrique sur $K_{0}[\mathfrak{a}]$.
On note alors \[\tilde{\beta}=v(Q_{n,i_{\beta}}(u_{1}(p),...,u_{n-1}(p),u_{n}(\beta)+a_{n,\beta}p^{\beta})),\]
\[\tilde{\varepsilon}=\max_{b\in\mathbb{N}}\left\lbrace\dfrac{\tilde{\beta}-\nu(\partial_{n,p^{b}}Q_{n,i_{\beta}})}{p^{b}}\right\rbrace\]
et\[\beta_{+}=\left \{ \begin{array}{ccc}  \min\lbrace\tilde\varepsilon,\varepsilon_{n,i_{\beta}}\rbrace & \textup{si} & \beta<\varepsilon_{n,i_\beta} \\  \min\lbrace\tilde\varepsilon,\varepsilon_{n,i_{\beta}+1}\rbrace  & \textup{si} & \beta=\varepsilon_{n,i_\beta} \end{array} \right.\]
Enfin, on pose: \[u_{n}(\beta_{+})=u_{n}(\beta)+a_{n,\beta}p^{\beta},\] 
ceci nous d\'efinit alors un nouveau d\'eveloppement de Puiseux partiel sur lequel on peut continuer la r\'ecurrence. En effet, remarquons que $\tilde\varepsilon\geqslant\beta$ car si $\beta<\varepsilon_{n,i_\beta}$, alors, par d\'efinition de $\nu_\beta$, $\tilde\beta\geqslant\nu_{\beta}(Q_{n,i_\beta})=\min\limits_{q\in\mathbb{N}}\left\lbrace \nu\left( \partial_{n,p^{q}}Q_{n,i_{\beta}}\right) +p^{q}\beta\right\rbrace$ et donc $\tilde\varepsilon\geqslant\beta$ par d\'efinition de $\tilde\varepsilon$; si $\beta=\varepsilon_{n,i_\beta}$, par le (3) du Lemme \ref{lemmebeta}, $\tilde\beta=\beta_{n,i_\beta}$ et donc, toujours par d\'efinition $\tilde\varepsilon\geqslant\beta$. Ainsi, le (1) de la D\'efinition \ref{defdevpart} est toujours v\'erifi\'e lorsque $\beta<\varepsilon_{n,i_\beta}$; si $\beta=\varepsilon_{n,i_\beta}$, c'est \'egalement vrai vu que, dans ce cas, $\tilde\beta=\beta_{n,i_\beta}$. Quant au (2), on vient de voir que $\tilde\beta\geqslant\nu_{\beta}(Q_{n,i_\beta})$ pour $\beta<\varepsilon_{n,i_\beta}$; si $\beta=\varepsilon_{n,i_\beta}$, comme la suite $(\beta_{n,i})_{i\in\Lambda_n}$ est croissante, on a, lorsque $\beta_+=\tilde\varepsilon$:
\[v(\iota_{\beta_+}(Q_{n,i_{\beta}+1}))\geqslant\tilde\beta=\beta_{n,i_\beta}.\]
On en d\'eduit que $v(\iota_{\beta_+}(Q_{n,i_{\beta}+1}))\geqslant\min\limits_{q\in\mathbb{N}}\left\lbrace \nu\left( \partial_{n,p^{q}}Q_{n,i_{\beta}+1}\right) +p^{q}\beta_+\right\rbrace$. Si $\beta_+=\varepsilon_{n,i_\beta+1}$, alors: \[v(\iota_{\beta_+}(Q_{n,i_{\beta}+1}))=\min\limits_{q\in\mathbb{N}}\left\lbrace \nu\left( \partial_{n,p^{q}}Q_{n,i_{\beta}+1}\right) +p^{q}\beta_+\right\rbrace\] par le (3) du Lemme \ref{lemmebeta}.
\end{enumerate}

\indent Pour achever notre r\'ecurrence transfinie, il nous faut consid\'erer le cas limite. Soient $\mathcal{W}$ un sous-ensemble bien ordonn\'e de $\Gamma_{1}$ n'ayant pas d'\'el\'ement maximal et $\lbrace a_{n,\gamma}\:\vert\:\gamma\in \mathcal{W}\rbrace$ tels que $\forall\:\beta\in \mathcal{W}$, $u_{n}(\beta)=\sum\limits_{\underset{\gamma < \beta}{\gamma\in \mathcal{W}}}a_{n,\gamma}p^{\gamma}$ soit un d\'eveloppement de Puiseux partiel de $u_{n}$. Notons $u_{n}(\mathcal{W})=\sum\limits_{\gamma\in \mathcal{W}}a_{n,\gamma}p^{\gamma}$.
\\Supposons d'abord que:
\[ \forall\: i\in\Lambda_{n},\:\exists\:\beta\in \mathcal{W},\:\varepsilon_{n,i}<\beta. \]
Alors, pour tout $i\in\Lambda_{n},\:i<i_{\beta}$ et donc l'ensemble $\left\lbrace Q_{n,i} \right\rbrace _{i\in\Lambda_{n}}$ forme un syst\`eme complet de polyn\^omes-cl\'es pour l'extension $K_{n-1}\hookrightarrow K_{n-1}(u_{n})$. Ainsi:
\[\forall\:f\in K_{n-1}[u_{n}],\:\exists \:i\in\Lambda_{n},\:\nu(f)=\nu_{n,i}(f).\]
On en d\'eduit donc, \`a l'aide de la d\'efinition du d\'eveloppement de Puiseux partiel, que:
\[\forall\:f\in K_{n-1}[u_{n}],\:\nu(f)=v(f(u_{1}(p),...,u_{n-1}(p),u_{n}(\mathcal{W}))).\]
Or tout $f\in R$ tel que $\nu(f)\in\Gamma_{1}$ s'\'ecrit $f=f'+f''$, avec $f'\in K_{n-1}[u_{n}]$ et $\nu_{0}(f'')>\nu(f)$.
\\On a alors:
\begin{align*}
v(f''(u_{1}(p),...,u_{n-1}(p),u_{n}(\mathcal{W})))>\nu(f)&=\nu(f')\\&=v(f'(u_{1}(p),...,u_{n-1}(p),u_{n}(\mathcal{W})))\\&=v(f(u_{1}(p),...,u_{n-1}(p),u_{n}(\mathcal{W}))).
\end{align*}
D'o\`u, pour tout $f\in R$ tel que $\nu(f)\in\Gamma_{1}$, on a:
\[\nu(f)=v(f(u_{1}(p),...,u_{n-1}(p),u_{n}(\mathcal{W}))).\]
Enfin, le m\^eme r\'esultat est vrai pour tout $f\in R\:\otimes_{R_{n-1}} K_{n-1}$ tel que $\nu(f)\in\Gamma_{1}$.
\\S'il existe un $f\in R$ tel que $\nu(f)\notin\Gamma_{1}$, alors l'ensemble $\Lambda_{n}$ contient un \'el\'ement maximal $\lambda$ et donc il existe un $\beta\in \mathcal{W}$ tel que $\varepsilon_{n,\lambda}<\beta$. Alors $f$ s'\'ecrit de mani\`ere unique sous la forme $f=Q_{n,\lambda}^{a}\tilde{f}$, o\`u $a\in\mathbb{N}$, $\tilde{f}\in R\:\otimes_{R_{n-1}} K_{n-1}$ tel que $\nu(\tilde{f})\in\Gamma_{1}$. Par le cas pr\'ec\'edent, $\nu(\tilde{f})=v(\tilde{f}(u_{1}(p),...,u_{n-1}(p),u_{n}(\mathcal{W})))$ et donc $\nu(f)=v(f(u_{1}(p),...,u_{n-1}(p),u_{n}(\mathcal{W})))$.
\\On d\'efinit alors $u_{n}(p)=u_{n}(\mathcal{W})$ et la construction du d\'eveloppement de Puiseux s'arr\^ete.
\\Supposons, pour terminer, que:
\[ \exists\: i\in\Lambda_{n},\:\forall\:\beta\in \mathcal{W},\:\varepsilon_{n,i}>\beta. \] 
On note alors\[i_{\mathcal{W}}=\min\lbrace i\in\Lambda_{n}\:\vert\:\forall\:\beta\in \mathcal{W},\:\varepsilon_{n,i}>\beta\rbrace,\] 
\[\tilde{\beta}=v(Q_{n,i_{\mathcal{W}}}(u_{1}(p),...,u_{n-1}(p),u_{n}(\mathcal{W}))),\] 

\[\tilde{\varepsilon}=\max_{b\in\mathbb{N}} \left\lbrace\dfrac{\tilde{\beta}-\nu(\partial_{n,p^{b}}Q_{n,i_{\mathcal{W}}})}{p^{b}}\right\rbrace\]
et\[\beta_{+}=\min\lbrace\tilde\varepsilon,\varepsilon_{n,i_{\mathcal{W}}}\rbrace.\]
Enfin, on pose $u_{n}(\beta_{+})=u_{n}(\mathcal{W})$, ceci nous d\'efinit alors un nouveau d\'eveloppement de Puiseux partiel sur lequel on peut continuer la r\'ecurrence.
\\ \qed

\section{Des r\'esultats de d\'ependance int\'egrale}
Les r\'esultats de cette section sont donn\'es dans le cas mixte, en rempla\c{c}ant $W$ par $k$, $\overline{W}$ par $\overline{k_{\nu}}$, $p$ par $t$ et $n$ par $n+1$, on obtient les m\^emes r\'esultats dans le cas \'equicaract\'eristique.
\begin{prop}
Soient $i\in\Lambda_{n}$, $\beta=\varepsilon_{n,i}$ (c'est-\`a-dire $i=i_{\beta}$), et $h\in R$. Alors:
\[\nu_{n,i}(h)=\min_{\alpha\in\mathbb{N}}\left\lbrace \nu (\partial_{n,\alpha} h)+\alpha\beta\right\rbrace =\min_{\alpha\in\mathbb{N}}\left\lbrace \nu_{n,i}(\partial_{n,\alpha}h)+\alpha\beta\right\rbrace.\]\label{min}
\end{prop}
\noindent\textit{Preuve}: Soient $h\in R$ et $\alpha\in\mathbb{N}$, par la Proposition \ref{prop10.1}, on a:
\[\nu_{n,i}(h)-\nu_{n,i}(\partial_{n,\alpha}h)\leqslant\alpha\beta.\]
On obtient alors:
\[\nu_{n,i}(h)\leqslant\min_{\alpha\in\mathbb{N}}\left\lbrace \nu_{n,i}(\partial_{n,\alpha}h)+\alpha\beta\right\rbrace\leqslant\min_{\alpha\in\mathbb{N}}\left\lbrace \nu (\partial_{n,\alpha} h)+\alpha\beta\right\rbrace.\]
Montrons que ces in\'egalit\'es sont des \'egalit\'es. Par la Remarque \ref{rmqval} (3), on a:
\[\nu_{n,i}(h)=\mu_{\beta}(h)=\nu_{\beta}(\iota_{\beta,T}(h)).\]
Soit $h=\sum\limits_{j=0}^{s_{n,i}}d_{n,j,i}Q_{n,i}^{j}$ le d\'eveloppement $(n,i)$-standard de $h$, on pose:
\[S=\lbrace j\in\lbrace 0,...,s_{n,i}\rbrace\:\vert\: j\beta + \nu(d_{n,j,i})=\nu_{n,i}(h)\rbrace.\]
Tout \'el\'ement $j\in S$ s'\'ecrit de la forme $j=p^{e}u$ o\`u $p$ ne divise pas $u$. Prenons un $j\in S$ tel que $p^{e+1}$ divise $j'$, pour tout $j'\in S$ tel que $j'<j$. On pose alors $\alpha=p^{b_{n,i}}j$ et \`a l'aide de la Proposition \ref{coro10.15} on en d\'eduit que:
\[\nu_{n,i}(h)-\nu(\partial_{n,\alpha}h)=\nu_{n,i}(h)-\nu_{n,i}(\partial_{n,\alpha}h)=\alpha\beta.\]
\qed
\\ \\ \indent Notons $\mathcal{A}$ le sous-anneau de $\overline{W}\left[ \left[ p^{\Gamma'}\right] \right] $ engendr\'e par $\iota (p),\iota (u_{1}),...,\iota (u_{n})$ et toutes leurs troncatures. Pour tout $j\in \lbrace r,...,n\rbrace$ et $\beta\in\Gamma\cup\lbrace\infty\rbrace$, on note $\mathcal{A}_{j,\beta}$ le sous-anneau de $\mathcal{A}$ engendr\'e par toutes les troncatures ouvertes de la forme $u_{j'}(\beta')$, o\`u $(j',\beta')<_{lex}(j,\beta)$ pour l'ordre lexicographique.
\begin{prop}
Soient $j\in \lbrace r,...,n\rbrace$, $\beta\in\Gamma\cup\lbrace\infty\rbrace$, $g,h\in\mathcal{A}_{j,\beta}[u_{j}(\beta)]$ et $\lambda\in\Gamma_{1}$. On suppose que $v(gh)<\lambda$.
\\Il existe alors $l\in\mathbb N$, $\lambda_{0}<...<\lambda_{l}$ et $\delta_{1}>...>\delta_{l}$ \'el\'ements de $\Gamma_{1}$ tels que:
\[(gh)(\lambda)=\sum\limits_{i=1}^{l}g[\lambda_{i-1},\lambda_{i}[\:h(\delta_{i}).\]
De plus, on peut choisir les suites $(\lambda_{i})_{0\leqslant i\leqslant l}$ et $(\delta_{i})_{1\leqslant i \leqslant l}$ de telle sorte que $\lambda_{l}\leqslant \lambda - v(h)$ et $\delta_{1}\leqslant \lambda -v(g)$.
\label{caltron}\end{prop}
\noindent\textit{Preuve}: Notons $supp(g)$ (resp. $supp(h)$) l'ensemble de tous les $\varepsilon\in\Gamma_{1}$ tels que le coefficient devant $p^{\varepsilon}$ de $g$ (resp. de $h$) soit non-nul. On va construire les deux suites cherch\'ees par r\'ecurrence.
\\On pose $\lambda_{0}=v(g)$ et $\delta_{1}=\lambda-\lambda_{0}$, par hypoth\`eses on a bien $\lambda_{0}\leqslant \lambda - v(h)$. Supposons maintenant que, pour tout $q\geqslant 1$, on ait construit $\lambda_{0}<...<\lambda_{q}$ et $\delta_{1}>...>\delta_{q}$ avec $\lambda_{q}\leqslant \lambda - v(h)$ et $\delta_{i}=\lambda-\lambda_{i-1}$, pour $1\leqslant i \leqslant q$. Posons alors:
\[B_{q}=\lbrace\varepsilon\in supp(g)\:\vert\:\exists\:\theta\in supp(h),\: \theta + \lambda_{q}<\lambda\leqslant\theta + \epsilon\rbrace.\]
Si $B_{q}=\emptyset$, on pose $l=q$ et la r\'ecurrence s'arr\^ete (remarquons que ceci arrive lorsque $\lambda_{q}\geqslant \lambda - v(h)$). De plus, par construction, on a l'\'egalit\'e:
\[(gh)(\lambda)=\sum\limits_{i=1}^{l}g[\lambda_{i-1},\lambda_{i}[\:h(\delta_{i}).\]
\\Si $B_{q}\neq\emptyset$, on pose $\lambda_{q+1}=\min \lbrace\lambda - v(h),\:\min  B_{q}\rbrace$ et $\delta_{q+1}=\lambda-\lambda_{q}$. Par d\'efinition de $\lambda_{q+1}$ et de $\delta_{q+1}$ et par hypoth\`ese de r\'ecurrence on a bien que $\lambda_{q}<\lambda_{q+1}$ et $\delta_{q}>\delta_{q+1}$. De plus, on remarque que:
\[\lbrace \lambda-\lambda_{q+1},\:\lambda-\lambda_{q}\rbrace\cap supp(h)\neq\emptyset.\]
On obtient alors une suite strictement d\'ecroissante d'ensembles:
\[supp\left( h(\lambda-\lambda_{1})\right) \supsetneqq ... \supsetneqq supp\left( h(\lambda-\lambda_{q+1})\right),\]
o\`u $supp\left( h(\lambda-\lambda_{q+1})\right)$ est un segment initial de $supp\left( h(\lambda-\lambda_{q})\right)$. Le processus s'arr\^ete donc au bout d'un nombre fini d'it\'erations, ceci entra\^inant la finitude des suites $(\lambda_{i})_{i}$ et $(\delta_{i})_{i}$.\\ \qed
\begin{coro}
Soient $j\in \lbrace r,...,n\rbrace$, $\beta\in\Gamma\cup\lbrace\infty\rbrace$, $g_{1},...,g_{s}\in\mathcal{A}_{j,\beta}[u_{j}(\beta)]$ et $\lambda\in\Gamma_{1}$. On suppose que $v(g_{1}...g_{s})<\lambda$. Alors:
\[(g_{1}...g_{s})(\lambda)=\sum\limits_{(i_{1},...,i_{s})\in N}\prod\limits_{j=1}^{s}g_{j}(\lambda_{i_{j}}^{(j)}),\]
o\`u $N\subset(\mathbb{N}^{*})^{s}$ est un ensemble fini et $\lambda_{i_{j}}^{(j)}\in\Gamma_{1\:+}$ sont tels que $\lambda_{i_{j}}^{(j)}\leqslant\lambda$ avec in\'egalit\'e stricte s'il existe $j'\in\lbrace 1,...,s\rbrace\backslash\lbrace j\rbrace$ tel que $v(g_{j'})>0$.
\label{prodfini}\end{coro}
\noindent\textit{Preuve}: Par r\'ecurence sur $s$ en appliquant la Proposition \ref{caltron}.\\ \qed
\begin{coro}
Soient $j\in \lbrace r,...,n\rbrace$, $\beta\in\Gamma\cup\lbrace\infty\rbrace$, $g_{1},...,g_{s}\in\mathcal{A}_{j,\beta}[u_{j}(\beta)]$ et $\lambda\in\Gamma_{1}$. On suppose que $v(g_{1}...g_{s})<\lambda$. Alors:
\[(g_{1}...g_{s})(\lambda)\in\mathcal{A}_{j,\beta}[u_{j}(\beta)].\]
\label{stabtron}\end{coro}
\noindent\textit{Preuve}: Par le Corollaire \ref{prodfini}, il suffit de montrer le r\'esultat pour $s=1$. Notons $g=g_{1}$ et montrons par r\'ecurrence sur $j\in\lbrace r,...,n\rbrace$ que, si $g\in\mathcal{A}_{j,\beta}[u_{j}(\beta)]$, alors $g(\lambda)\in\mathcal{A}_{j,\beta}[u_{j}(\beta)]$, pour $\beta\in\Gamma\cup\lbrace\infty\rbrace$ et $\lambda\in\Gamma_{1}$ fix\'es.
\\Pour $j=r$ on a:
\[ \mathcal{A}_{r,\beta}[u_{r}(\beta)]=\left \{ \begin{array}{ccc}  \overline{W}\left[ p^{\nu(p)},p^{\nu(u_{1})},...,p^{\nu(u_{r})}\right]  & \textup{si} & \nu(u_{r})<\beta\\  \overline{W}\left[ p^{\nu(p)},p^{\nu(u_{1})},...,p^{\nu(u_{r-1})}\right]   & \textup{si} & \nu(u_{r})\geqslant\beta \end{array} \right.\]
et donc si $g\in\mathcal{A}_{r,\beta}[u_{r}(\beta)]$ et $v(g)<\lambda$ alors $g(\lambda)\in\mathcal{A}_{r,\beta}[u_{r}(\beta)]$.
\\Supposons $j>r$ et le r\'esultat vrai pour $j-1$. Soit $g\in\mathcal{A}_{j,\beta}[u_{j}(\beta)]$, on peut \'ecrire $g$ comme un polyn\^ome en $u_{j}(\beta)$ \`a coefficients dans $\mathcal{A}_{j,\infty}$. On applique alors le Corollaire \ref{prodfini} \`a chaque mon\^ome de $g(\lambda)$. Par hypoth\`ese de r\'ecurrence, toutes les troncatures ouvertes des coefficients de $g$ sont dans $\mathcal{A}_{j,\infty}$, ainsi $g(\lambda)\in\mathcal{A}_{j,\beta}[u_{j}(\beta)]$.\\ \qed
\\ \\ \indent Pour $j\in \lbrace r,...,n\rbrace$ et $\beta\in\Gamma\cup\lbrace\infty\rbrace$, consid\'erons les morphismes d'anneaux $\tau_{j,\beta},\iota_{j,\beta}:R_{j}\rightarrow\overline{W}\left[ \left[ p^{\Gamma'}\right] \right] $ d\'efinis par:
\begin{align*}
\tau_{j,\beta}(p)&=\iota_{j,\beta}(p)=\iota(p); \\
\tau_{j,\beta}(u_{i})&=\iota_{j,\beta}(u_{i})=\iota(u_{i}),\:\forall\:i\in\lbrace 1,...,j-1\rbrace ; \\
\iota_{j,\beta}(u_{j})&=u_{j}(\beta);\\
 \tau_{j,\beta}(u_{j})&= u_{j}[\beta].
\end{align*}
Remarquons que $\iota_{n,\beta}=\iota_{\beta}$. Rappelons la notation $R_j=W\left[\left[u_1,...,u_j\right]\right]$.
\begin{prop}
Soient $j\in \lbrace r,...,n\rbrace$, $\beta\in\Gamma\cup\lbrace\infty\rbrace$, $f\in R_{j}\backslash R_{j-1}$ et $\lambda\in\Gamma_{1}$. Alors:
\[ \iota_{j,\beta}(f) (\lambda)\in\mathcal{A}_{j,\beta}[u_{j}(\beta)].\]
\label{stab}\end{prop}
\noindent\textit{Preuve}: Remarquons tout d'abord que l'on peut remplacer $f$ par une de ses approximations $(p,u_{1},...,u_{j})$-adiques choisies dans $W\left[ u_{1},...,u_{j}\right] $ de telle sorte que l'on ne modifie pas $ \iota_{j,\beta}(f) (\lambda)$. Soit $f\in W\left[ u_{1},...,u_{j}\right]$ choisi ainsi, on a alors:
\[ \iota_{j,\beta}(f) \in W \left[ \iota(p),\iota(u_{1}),...,\iota(u_{j-1}), u_{j}(\beta) \right] \subset\mathcal{A}_{j,\beta}[u_{j}(\beta)].\]
On conclut en appliquant le Corollaire \ref{stabtron}. \\ \qed
\\ \\ \indent Dans ce qui suit nous allons donner une description explicite de $ \tau_{j,\beta}(f) (\lambda)$ et $ \iota_{j,\beta}(f) (\lambda)$ pour un $\lambda$ que nous pr\'eciserons par la suite.
\\Soient $j\in \lbrace r,...,n\rbrace$, $\beta\in\Gamma_{1}\cup\lbrace\infty\rbrace$ et $f\in R_{j}\backslash R_{j-1}$. Notons \[i_{j,\beta}=\min\lbrace i\in\Lambda_{j}\:\vert\:\beta\leqslant\varepsilon_{j,i}\rbrace\]
et par convention, si $\lbrace i\in\Lambda_{j}\:\vert\:\beta\leqslant\varepsilon_{j,i}\rbrace=\emptyset$, on prendra $i_{j,\beta}=\Lambda_{j}$ (c'est-\`a-dire le plus petit ordinal strictement plus grand que n'importe quel \'el\'ement de $\Lambda_{j}$). Remarquons que $i_{n,\beta}=i_{\beta}$. On pose alors:
\[\lambda(f,\beta)=\min\lbrace \nu_{j,i_{j,\beta}}(\partial_{j,b}f)+b\beta\:\vert\:b\in\mathbb{N}^{*}\rbrace;\]
\[U=\lbrace b\in\mathbb{N}^{*}\:\vert\:\nu_{j,i_{j,\beta}}(\partial_{j,b}f)+b\beta=\lambda(f,\beta)\rbrace.\]
\begin{rem}
\textup{Comme $R_{j}$ est noeth\'erien, on a le fait suivant:
\[\exists\:\overline{b}\in\mathbb{N}^{*},\:\partial_{j,b}f\in\left( \partial_{j,0}f,...,\partial_{j,\overline{b}}f\right),\:\forall\:b>\overline{b}.\]
Ainsi $U\subset\lbrace 0,...,\overline{b}\rbrace$ est un ensemble fini.}
\end{rem}
Par abus de notation, on notera $in_{v}(\partial_{j,b}f)$ le mon\^ome de plus petit degr\'e de $\iota(\partial_{j,b}f)$ dans $\overline{W}\left[ \left[ p^{\Gamma'}\right] \right] $. On appelle $U_{0}$ l'ensemble des $b\in U$ tels que $in_{\varepsilon_{j,i_{j,\beta}}}(T)$ n'appara\^it pas dans $in_{\varepsilon_{j,i_{j,\beta}}}\left( \iota_{\varepsilon_{j,i_{j,\beta}},T}(\partial_{j,b}f)\right) $. En rempla\c{c}ant $n$ par $j$ dans le Lemme \ref{rmqval} (4), on a, pour tout $b\in U_{0}$:
\[\nu_{j,i_{j,\beta}}(\partial_{j,b}f)=\nu(\partial_{j,b}f)=v\left( \iota_{\varepsilon_{j,i_{j,\beta}}}(\partial_{j,b}f)\right);\]
\[in_{v}(\partial_{j,b}f)=in_{v}\left( \iota_{\varepsilon_{j,i_{j,\beta}}}(\partial_{j,b}f)\right).\]
\begin{rem}
\textup{Soit $b\in U$, on a alors:
\[b\in U_{0}\;\textup{ssi}\;\exists\:i_{0}<i_{j,\beta},\:\nu_{j,i_{0}}(\partial_{j,b}f)=\nu(\partial_{j,b}f).\]
Comme $U$ et $U_{0}$ sont des ensembles finis, le m\^eme $i_{0}$ peut \^etre choisi tel que $\nu_{j,i_{0}}(\partial_{j,b}f)=\nu(\partial_{j,b}f)$ et ceci pour tout $b\in U_{0} $. }
\end{rem}
Pour $i_{j,\beta}$ ordinal limite, prenons $i_{0}\in\Lambda_{j}$ satisfaisant les conditions suivantes:
\begin{enumerate}
\item $i_{0}<i_{j,\beta}$;
\item $\forall\:b\in U_{0}$, $\nu_{j,i_{0}}(\partial_{j,b}f)=\nu(\partial_{j,b}f)$;
\item $\forall\: i\in\Lambda_{j}$, $i_{0}<i<i_{j,\beta}$, $\forall\: b\in U_{0}$,
\[v(\iota(\partial_{j,b}f)-in_{v}(\partial_{j,b}f))-\nu(\partial_{j,b}f)>\beta-\varepsilon_{j,i}\:;\]
\item $\forall\:i\in\Lambda_{j}$, $i_{0}<i<i_{j,\beta}$, $\forall\:b\in\lbrace 0,...,\overline{b}\rbrace\setminus U_{0}$,
\[\nu(\partial_{j,b}f)+b\varepsilon_{j,i}>\lambda(f,\beta).\]
\end{enumerate}
Enfin, notons:
\[\Delta u_{j}=u_{j}[\beta]-u_{j}[\varepsilon_{j,i_{0}}];\]
\[\Delta u_{j}(\beta)=u_{j}(\beta)-u_{j}[\varepsilon_{j,i_{0}}].\]
\begin{prop}
Il existe deux polyn\^omes $F_{\beta},\tilde{F}_{\beta}\in\mathcal{A}_{j,\beta}\left[ X\right] $ de la forme:
\[F_{\beta}\left( X\right) =F_{0}+\sum\limits_{b\in U_{0}}in_{v}(\partial_{j,b}f)X^{b};\]
\[\tilde{F}_{\beta}\left( X\right) =\tilde{F}_{0}+\sum\limits_{b\in U_{0}}in_{v}(\partial_{j,b}f)X^{b};\]
tels que:
\begin{enumerate}
\item $F_{0},\tilde{F}_{0}\in\mathcal{A}_{j,\beta}$;
\item $\tau_{j,\beta}(f)\left[ \lambda(f,\beta)\right] =F_{\beta}\left( \Delta u_{j}\right)$;
\item $\iota_{j,\beta}(f)\left( \lambda(f,\beta)\right) =\tilde{F}_{\beta}\left( \Delta u_{j}(\beta)\right)$.
\end{enumerate}\label{taylor}
\end{prop}
\noindent\textit{Preuve}: Dans ce qui suit, on adoptera la convention $\varepsilon_{j,\Lambda_{j}}=\infty$, pour tout $j\in\lbrace r+1,...,n\rbrace$. Supposons d'abord que $\beta\leqslant\varepsilon_{j,1}$, alors, $u_{j}[\beta]=\Delta u_{j}=0$ et on pose $F_{0}=0$. Dans ce cas, (2) est trivialement montr\'e et on proc\`ede de la m\^eme mani\`ere pour montrer (3).
\\Supposons que $\beta>\varepsilon_{j,1}$. Soit $\varepsilon\in\Gamma_{1}\otimes_{\mathbb{Z}}\mathbb{Q}$, $\varepsilon>0$ suffisamment petit, notons $\partial_{j,b}f\left( u_{j}[\varepsilon_{j,i_{0}}]\right) = \iota_{j,\varepsilon_{j,i_{0}}+\varepsilon}(\partial_{j,b}f)$ et prenons le d\'eveloppement de Taylor de $\iota(f)$ en $u_{j}[\beta]$ en la $j$-\`eme variable. On peut alors \'ecrire:
\[\tau_{j,\beta}(f)=\sum\limits_{b\in\mathbb{N}}\partial_{j,b}f\left( u_{j}[\varepsilon_{j,i_{0}}]\right)\left( \Delta u_{j}\right) ^{b}.\]
Par les points $(3)$ et $(4)$ dans le choix de $i_{0}$, les termes de la forme:
\[\left( \iota(\partial_{j,b}f)-in_{v}(\partial_{j,b}f)\right) \left( u_{j}[\varepsilon_{j,i_{0}}]\right)\left( \Delta u_{j}\right) ^{b},\; \textup{pour}\;b\in U_{0};\]
\[\textup{et}\;\partial_{j,b}f\left( u_{j}[\varepsilon_{j,i_{0}}]\right)\left( \Delta u_{j}\right) ^{b},\; \textup{pour}\;b\notin U_{0};\]
n'interviennent pas dans $\tau_{j,\beta}(f)\left[ \lambda(f,\beta)\right]$. Ainsi, \`a l'aide de la Proposition \ref{stab}, on a:
\[\tau_{j,\beta}(f)\left[ \lambda(f,\beta)\right]-\sum\limits_{b\in U_{0}}in_{v}(\partial_{j,b}f)\left( \Delta u_{j}\right)^{b}=\iota_{j,\varepsilon_{j,i_{0}}+\varepsilon}(f)\in\mathcal{A}_{j,\varepsilon_{j,i_{0}}+\varepsilon}\left[ u_{j}\left( \varepsilon_{j,i_{0}}+\varepsilon\right) \right] .\]
On en d\'eduit donc que $F_{0}:=\tau_{j,\beta}(f)\left[ \lambda(f,\beta)\right]-\sum\limits_{b\in U_{0}}in_{v}(\partial_{j,b}f)\left( \Delta u_{j}\right)^{b}\in\mathcal{A}_{j,\beta}$.
De la m\^eme mani\`ere on montre, en faisant un d\'eveloppement de Taylor de $\iota(f)$ en $u_{j}(\beta)$, que
$\tilde{F}_{0}:=\iota_{j,\beta}(f)\left( \lambda(f,\beta)\right)-\sum\limits_{b\in U_{0}}in_{v}(\partial_{j,b}f)\left( \Delta u_{j}(\beta)\right)^{b}\in\mathcal{A}_{j,\beta}$. Ainsi, $(1)$, $(2)$, et $(3)$ sont d\'emontr\'es. Pour conclure il suffit de montrer que:
\[\forall\:b\in U_{0},\:in_{v}(\partial_{j,b}f)\in\mathcal{A}_{j,\beta}.\]
Soit $b\in U_{0}$ et notons $g=\partial_{j,b}f$. On pose $i_{0}(b)=\min\lbrace i_{0}\in\Lambda_{j}\:\vert\:\nu_{j,i_{0}}(g)=\nu(g)\rbrace$. Par la Proposition \ref{coro10.15}, on a:
\[\nu_{j,i_{0}(b)}(g)=\nu(g)\leqslant\nu(\partial_{j,q}g)+q\varepsilon_{j,i_{0}(b)},\:\forall\:q\in\mathbb{N}.\]
Par minimalit\'e de $i_{0}(b)$, il existe $q>0$ tel que l'on ait \'egalit\'e dans l'in\'egalit\'e pr\'ec\'edente; prenons alors le plus petit $q$ de la sorte. Le choix d'un tel $q$ entra\^ine que $\nu(\partial_{j,q}g)=\nu_{j,i_{0}(b)}(\partial_{j,q}g)$, en effet, sinon on aurait:
\[\nu_{j,i_{0}(b)}(g)\leqslant\nu_{j,i_{0}(b)}(\partial_{j,q}g)+q\varepsilon_{j,i_{0}(b)}<\nu(\partial_{j,q}g)+q\varepsilon_{j,i_{0}(b)}=\nu(g),\]
ce qui contredit le choix de $q$.
Pour $\varepsilon\in\Gamma_{1}$, $\varepsilon>0$ suffisamment petit, on a:
\[\lambda\left( g,\varepsilon_{j,i_{0}(b)}+\varepsilon\right) =\nu(\partial_{j,q}g)+q\left( \varepsilon_{j,i_{0}(b)}+\varepsilon\right) .\]
Fixons nous alors $\varepsilon\in\Gamma_{1}$, $\varepsilon>0$ suffisamment petit tel que $in_{v}(\iota(g))=\iota(g)\left( \lambda\left( g,\varepsilon_{j,i_{0}(b)}+\varepsilon\right)\right) $ et $\varepsilon_{j,i_{0}(b)}+\varepsilon<\min\lbrace\beta,\varepsilon_{j,i_{0}(b)+1}\rbrace$. Par la Proposition \ref{stab}, on en d\'eduit que:
\begin{align*}
in_{v}(g)&=\iota(g)\left( \lambda\left( g,\varepsilon_{j,i_{0}(b)}+\varepsilon\right)\right)\\&=\iota_{j,\varepsilon_{j,i_{0}(b)}+\varepsilon}(g)\left( \lambda\left( g,\varepsilon_{j,i_{0}(b)}+\varepsilon\right)\right)\\&\in\mathcal{A}_{j,\varepsilon_{j,i_{0}(b)}+\varepsilon}\left[ u_{j}\left( \lambda\left( g,\varepsilon_{j,i_{0}(b)}+\varepsilon\right)\right) \right]\subset\mathcal{A}_{j,\beta}.
\end{align*}\qed
\begin{prop}
Soient $j\in \lbrace r,...,n\rbrace$, $\beta\in\Gamma_{1}\cup\lbrace\infty\rbrace$. Si $i_{j,\beta}<\Lambda_{j}$, alors $u_{j}(\beta)$ est entier sur $\mathcal{A}_{j,\beta}$ et une relation de d\'ependance int\'egrale est donn\'ee par:
\[\tilde{F}_{0}\left( \lambda\left( Q_{j,i_{j,\beta}},\beta\right) \right) +\sum\limits_{b\in U_{0}}in_{v}\left( \partial_{j,b}Q_{j,i_{j,\beta}}\right) \left( \Delta u_{j}(\beta)\right) ^{b}=0.\]
\label{ent}\end{prop}
\noindent\textit{Preuve}: Soient $j\in \lbrace r,...,n\rbrace$, $\beta\in\Gamma_{1}\cup\lbrace\infty\rbrace$ et supposons que $i_{j,\beta}<\Lambda_{j}$. Par la Proposition \ref{taylor} et par construction des d\'eveloppements de Puiseux, on a:
\[\tilde{F}_{\beta}\left( \Delta u_{j}(\beta)\right) =\iota_{j,\beta}\left( Q_{j,i_{j,\beta}}\right) \left( \lambda\left( Q_{j,i_{j,\beta}},\beta\right)\right) =0.\]
De plus, $\max U_{0}=d_{u_{j}}^{\circ}\left( Q_{j,i_{j,\beta}}\right)$ donc la relation pr\'ec\'edente est bien une relation de d\'ependance int\'egrale.
\qed
\\ \\ \indent Pour $j\in\lbrace r+1,...,n\rbrace$, notons $\mathcal{A}_{j}$ le sous-anneau de $\mathcal{A}$ engendr\'e par $\mathcal{A}_{j-1,\infty}$ et tous les \'el\'ements de la forme $u_{j}(\varepsilon_{j,i})$, $i\in\Lambda_{j}$.
\begin{coro}
Soit $j\in\lbrace r+1,...,n\rbrace$, pour tout $\beta\in\Gamma_{1}$ tel que $i_{j,\beta}<\Lambda_{j}$ on a:
\begin{enumerate}
\item $\mathcal{A}_{j,\beta}$ est une extension enti\`ere de $W\left[ \iota(p),\iota(u_{1}),...,\iota(u_{j-1})\right]$;
\item $\mathcal{A}_{j}$ est une extension enti\`ere de $W\left[ \iota(p),\iota(u_{1}),...,\iota(u_{j-1})\right]$;
\item $\mathcal{A}$ est une extension enti\`ere de $W\left[ \iota(p),\iota(u_{1}),...,\iota(u_{n})\right]$.
\end{enumerate}\label{extent}
\end{coro}
\noindent\textit{Preuve}: La premi\`ere assertion d\'ecoule de la construction des d\'eveloppements de Puiseux, la deuxi\`eme est une appliction imm\'ediate de la Proposition \ref{ent}, enfin, la troisi\`eme provient des deux pr\'ec\'edentes.\\ \qed
%\begin{prop}
%Soit $j\in\lbrace r+1,...,n\rbrace$, alors $\iota(u_{j})$ est transcendant sur $W\left[ \iota(p),\iota(u_{1}),...,\iota(u_{j-1})\right]$.
%\label{tran}\end{prop}
%\noindent\textit{Preuve}: Soit $P(X)\in W\left[ \iota(p),\iota(u_{1}),...,\iota(u_{j-1})\right]\left[ X\right] \setminus\lbrace 0\rbrace$, on peut construire $f\in W\left[ u_{1},...,u_{j-1}\right]\setminus\lbrace 0\rbrace $ tel que $\iota(f)=P\left(\iota( u_{j})\right) $. Or $f\neq 0$ donc:
%\[\nu(f)=v(\iota(f))\in\Gamma.\]
%Ainsi, $P(\iota(u_{j}))=\iota(f)\neq 0$.\\ \qed
\begin{coro}
Pour  $j\in\lbrace r+1,...,n\rbrace$, $\iota(u_{j})$ est transcendant sur $\mathcal{A}_{j}$.
\end{coro}
\noindent\textit{Preuve}: Soit $P(X)\in\mathcal{A}_{j}\left[ X\right]\setminus\lbrace 0\rbrace $. Par le Corollaire \ref{extent} (2), on sait que $\mathcal{A}_{j}$ est une extension enti\`ere de $W\left[ \iota(p),\iota(u_{1}),...,\iota(u_{j-1})\right]$, il existe donc un $W\left[ \iota(p),\iota(u_{1}),...,\iota(u_{j-1})\right]$-module de type fini contenu dans $\mathcal{A}_{j}$ et contenant tous les coefficients de $P$. On peut donc \'ecrire:
\[P(X)=\sum\limits_{b=1}^{c}P_{b}(X)a_{j,b}\]
o\`u $a_{j,b}\in\mathcal{A}_{j}$ et $P_{b}(X)\in W\left[ \iota(p),\iota(u_{1}),...,\iota(u_{j-1})\right][X]$, pour $b\in\lbrace 1,...,c\rbrace$. Comme $P\neq 0$, il existe $b_{0}\in\lbrace 1,...,c\rbrace$ tel que $P_{b_{0}}(X)\neq 0$. Comme $\iota(u_{j})$ est transcendant sur  $W\left[ \iota(p),\iota(u_{1}),...,\iota(u_{j-1})\right]$, alors $P_{b_{0}}\left( \iota(u_{j})\right) \neq 0$ et donc $P\left( \iota(u_{j})\right) \neq 0$.\\ \qed

%\nocite{*}
\bibliographystyle{cdraifplain}
\bibliography{biblio}
\end{document}